\documentclass{amsart}

\usepackage{amssymb,amsfonts, latexsym,amsmath,hhline,array,longtable}
\usepackage{pdfsync,color, comment,colortbl, mathrsfs,stmaryrd,cite,graphicx}

\usepackage{ifpdf}
\ifpdf \usepackage[colorlinks=true, citecolor=blue, linkcolor=blue, urlcolor=blue]{hyperref} \fi

\newcommand{\cal}{\mathcal}

\newtheorem{formula}{}[section]
\newtheorem{definition}[formula]{Definition}
\newtheorem{corollary}[formula]{Corollary}
\newtheorem{remark}[formula]{Remark}
\newtheorem{lemma}[formula]{Lemma}
\newtheorem{theorem}[formula]{Theorem}

\def\thrm{\begin{theorem}}
\def\thrml#1{\begin{theorem}\label{#1}}
\def\ethrm{\end{theorem}}
\def\rmrk{\begin{remark}}
\def\rmrkl#1{\begin{remark}\label{#1}}
\def\ermrk{\end{remark}}
\def\dfntn{\begin{definition}}
\def\dfntnl#1{\begin{definition}\label{#1}}
\def\edfntn{\end{definition}}
\def\nmrt{\begin{enumerate}}
\def\enmrt{\end{enumerate}}
\def\tm#1{\item[{\rm (#1)}]}
\def\qtnl#1{\begin{equation}\label{#1}}
\def\eqtn{\end{equation}}
\def\lmm{\begin{lemma}}
\def\lmml#1{\begin{lemma}\label{#1}}
\def\elmm{\end{lemma}}
\def\crllr{\begin{corollary}}
\def\crllrl#1{\begin{corollary}\label{#1}}
\def\ecrllr{\end{corollary}}
\def\css{\begin{cases}}
\def\ecss{\end{cases}}
\def\prf{\begin{proof}}
\def\eprf{\end{proof}}

\def\cG{{\cal G}}

\def\cX{{\cal X}}
\def\cY{{\cal Y}}

\def\mN{{\mathbb N}}

\def\fS{{\mathfrak S}}
\def\fX{{\mathfrak X}}
\def\fY{{\mathfrak Y}}

\DeclareMathOperator{\aut}{Aut}

\DeclareMathOperator{\diag}{Diag}

\DeclareMathOperator{\id}{id}
\DeclareMathOperator{\im}{im}

\DeclareMathOperator{\inv}{Inv}

\DeclareMathOperator{\iso}{Iso}

\DeclareMathOperator{\mon}{Mon}

\DeclareMathOperator{\pr}{pr}

\DeclareMathOperator{\rad}{Rad}

\DeclareMathOperator{\sym}{Sym}
\DeclareMathOperator{\WL}{WL}

\def\bone{{\bf 1}}
\def\grp#1{\langle {#1}\rangle}

\def\qaq{\quad\text{and}\quad}

\def\mltst#1{\{\hspace{-3pt}\{#1\}\hspace{-3pt}\}}
\def\ov{\overline}

\def\wh{\widehat}

\begin{document}

\title{On the WL-dimension of circulant graphs of prime power order}
\author{Ilia Ponomarenko}
\address{Steklov Institute of Mathematics at St. Petersburg, Russia}
\email{inp@pdmi.ras.ru}
\thanks{}
\date{}

\begin{abstract}
The $\WL$-dimension of a graph $X$ is the smallest positive integer~$m$ such that the $m$-dimensional Weisfeiler-Leman algorithm correctly tests the isomorphism between $X$ and any other graph. It is proved that the $\WL$-dimension of any circulant graph of prime power order is at most~$3$, and this bound cannot be reduced. The proof is based on using theories of coherent configurations and Cayley schemes over a cyclic group.
\end{abstract}

\maketitle

\section{Introduction}
A (not necessarily undirected) graph is said to be \emph{circulant} if it is a Cayley graph over a cyclic group. A main motivation of the present paper is the following computational problem: given a circulant graph $\cG$ and a graph $\cG'$ test whether $\cG$ is isomorphic to $\cG'$. To the best of the author's knowledge, there are two polynomial-time algorithms related to this problem, namely, \cite{EvdP2004a} and~\cite{Muz2004}. The first  solves the problem for all $\cG'$, but largely based on computational group theory, whereas the second  solves the problem only for circulant~$\cG'$ but is more or less  combinatorial. It is natural to ask if the problem can be solved in a purely combinatorial way for all~$\cG'$?

One of the most famous pure combinatorial methods for testing graph isomorphism was introduced  in~\cite{WLe68}. This method first uses the classical Weisfeiler-Leman algorithm, which, for a given graph~$\cG$, constructs a (in a sense) canonical coloring of all pairs of vertices. Then the resulting coloring is compared with similar coloring constructed for a graph~$\cG'$;  the graphs $\cG$ and $\cG'$ are declared isomorphic if and only if the colorings have the same sets of colors (this can be checked quite easily). This method tests  isomorphism of~$\cG$ and~$\cG'$ correctly in many (but not all) cases \cite{Pikhurko2011}.  

A generalization of the above method is obtained if the classical Weisfeiler-Leman algorithm is replaced by its $m$-dimensional ($m\ge 3$) analog, in which the $m$-tuples  of vertices are canonically colored. In this case, for every graph $\cG$ there is a minimal $m$ with the following property: any graph~$\cG'$ such that the set of colors in the canonical coloring of the $m$-tuples of the vertices of~$\cG'$ is equal to that for $\cG$, is isomorphic to~$\cG$. This minimal $m$ is called the \emph{Weisfeiler-Leman dimension} (the \emph{$\WL$-dimension}) of~$\cG$ and is denoted by $\dim_{\scriptscriptstyle\WL}(\cG)$. Now we can refine the question posed in the first paragraph as follows.\medskip

{\bf Question:}  {\it is it true that there exists $m\in\mN$ such that $\dim_{\scriptscriptstyle\WL}(\cG)\le m$ for every circulant graph~$\cG$.} 
\medskip

Although we cannot answer this question in full, our main result, presented by the theorem below, says that for circulant graphs of prime power order (i.e., those with the number of vertices equal to the power of a prime), as such a constant one can take~$3$.

\thrml{280222j}
The WL-dimension of every circulant graph of prime power order is at most~$3$.
\ethrm

Every Paley graph $\cG$ of prime order~$p$ is circulant. It is also known (see, e.g., \cite[Subsec.~4.5]{CP2019}) that for some $p$,  one can find a graph~$\cG'$ of order~$p$ such that $\cG$ and~$\cG'$ are nonisomorphic but have the same sets of colors in the colorings constructed by the classical Weisfeiler-Leman algorithm. It follows that $\dim_{\scriptscriptstyle\WL}(\cG)>2$, which shows that the estimate in Theorem~\ref{280222j} is sharp. An infinite family of such examples is obtained if one replaces $\cG$ and $\cG'$ by their lexicographic powers.  

Let us discuss some ideas underlying the proof of Theorem~\ref{280222j}. The classical Weisfeiler-Leman algorithm can be considered as an (efficiently computable) functor which given a graph~$\cG$  defines a coherent configuration $\cX=\WL(\cG)$  (the exact definitions are in Sections~\ref{160322a} and~\ref{120622s}). It was proved in~\cite{Fuhlbr2018a} that $\dim_{\scriptscriptstyle\WL}(\cG)\le 2$ if and only if the coherent configuration $\cX$ is separable, i.e., every algebraic isomorphism from $\cX$ to another coherent configuration is induced by isomorphism.

A weakening of the property of a coherent configuration to be separable is to consider not all algebraic isomorphisms, but only those that have one-point extensions. In this way, we arrive at a concept of {\it sesquiseparable} coherent configuration, which is introduced and studied in Section~\ref{160322f}. In fact, if $\cG$ is a vertex-transitive graph and the coherent configuration $\cX$ is sesquiseparable, then $\dim_{\scriptscriptstyle\WL}(\cG)\le 3$. We show that $\cX$ is sesquiseparable if $\cG$ is a circulant graph of prime power order.

Now let $\cG$ be a circulant graph. Then the coherent configuration $\cX=\WL(\cG)$ is a circulant scheme (a Cayley scheme over a cyclic group; for details, see Section~\ref{050322f}) and its structure is fairly well controlled by the Leung-Man theory (see, for example,~\cite{EvdP2003}). In particular, $\cX$ can be constructed from trivial and normal circulant schemes by using tensor and wedge products. The normal circulant schemes and the operation of wedge (or generalized wreath) product of (association) schemes have been introduced and studied in~\cite{EvdP2003} and~\cite{Muzychuk2009a}, respectively. It seems that a recent article~\cite{Herman2022} opens up a way to formulate the Leung-Man theory directly in terms of graph theory, bypassing association schemes.

When $\cG$ is of prime power order, the tensor product is irrelevant and we prove that $\cX$ is sesquiseparable by induction on the length of a decomposition of $\cX$ into a series of the wedge products (Lemma~\ref{010322a}). The central part of the proof  is focused on studying algebraic and combinatorial isomorphisms of the wedge product and establishing a sufficient condition for the wedge product to be sesquiseparable (see Section~\ref{090222e} and Remark~\ref{190322a}).

We conclude the introduction with a brief remark about Theorem~\ref{150222e} obtained in the course of the proof and which is of independent interest. In fact, this theorem says that Theorem~\ref{280222j} holds for much wider class of circulant graphs. Among these graphs are those $\cG$ for which the  scheme~$\cX$ is the wedge product of normal circulant schemes (such schemes cover the half of counterexamples constructed in~\cite{EvdP2002}). 

\section{Preliminaries}\label{160322a}

In order to make the paper as selfcontained as possible, we present in this section a necessary background of the theory of coherent configurations. Our notation is compatible with the notation in monograph~\cite{CP2019}; the proofs of the statements below can be found there.
 
 \subsection{Notation}
 Throughout the paper, $\Omega$ denotes a finite set. For $\Delta\subseteq \Omega$, the Cartesian product $\Delta\times\Delta$ and its diagonal are denoted by~$\bone_\Delta$ and $1_\Delta$, respectively. For a relation $s\subseteq\bone_\Omega$, we set $s^*=\{(\alpha,\beta): (\beta,\alpha)\in s\}$, $\alpha s=\{\beta\in\Omega:\ (\alpha,\beta)\in s\}$ for all $\alpha\in\Omega$, and define $\grp{s}$ as the minimal (with respect to inclusion) equivalence relation on~$\Omega$, containing~$s$. For any collection $S$ of relations, we denote by $S^\cup$ the set of all unions of elements of~$S$, and consider $S^\cup$ as a poset with respect to inclusion.
 
The set of classes of  an equivalence relation $e$ on~$\Omega$ is denoted by $\Omega/e$. For $\Delta\subseteq\Omega$, we set $\Delta/e=\Delta/e_\Delta$ where $e_\Delta=\bone_\Delta\cap e$. If the classes of $e_\Delta$ are singletons, $\Delta/e$ is identified with $\Delta$. Given a relation $s\subseteq \bone_\Omega$, we put
\qtnl{130622a}
s_{\Delta/e}=\{(\Gamma,\Gamma')\in \bone_{\Delta/E}: s^{}_{\Gamma,\Gamma'}\ne\varnothing\},
\eqtn
where $s^{}_{\Gamma,\Gamma'}=s\cap (\Gamma\times \Gamma')$. We also abbreviate $s^{}_\Gamma=s^{}_{\Gamma,\Gamma}$. Among all  equivalence relations $e$ on~$\Omega$, such that
\qtnl{kgu}
s=\bigcup_{(\Delta,\Delta')\in s_{\Delta/e}}\Delta\times \Delta',
\eqtn
there is  the largest (with respect to inclusion) one, which is denoted by  $\rad(s)$ and called the {\it radical} of~$s$. Obviously, $\rad(s)\subseteq\grp{s}$.

For a set $B$ of bijections $f:\Omega\to\Omega'$, subsets $\Delta\subseteq \Omega$ and $\Delta'\subseteq \Omega'$, equivalence relations $e$ and $e'$ on $\Omega$ and  $\Omega'$, respectively, we put
$$
B^{\Delta/e,\Delta'/e'}=\{f^{\Delta/e}:\ f\in B,\ \Delta^f=\Delta',\ e^f=e'\},
$$
where $f^{\Delta/e}$ is the bijection from $\Delta/e$ onto $\Delta'/e'$ induced by~$f$; we also  abbreviate $B^{\Delta/e}=B^{\Delta/e,\Delta'/e}$ if $\Delta'$ is clear from context.

\subsection{Coherent configurations}
Let $S$ be a partition of $\Omega^2$. A pair $\mathcal{X}=(\Omega,S)$ is called a \emph{coherent configuration} on $\Omega$ if 
\nmrt
\tm{C1}  $1_\Omega\in S^\cup$,
\tm{C2} $s^*\in S$ for all $s\in S$, 
\tm{C3} given $r,s,t\in S$, the number $c_{rs}^t=|\alpha r\cap \beta s^{*}|$ does not depend on $(\alpha,\beta)\in t$. 
\enmrt
The number $|\Omega|$ is called the {\it degree} of~$\cX$. We say that $\cX$ is {\it trivial} if  $S$ consists of~$1_\Omega$ and its complement (unless $\Omega$ is not a singleton), {\it homogeneous} or a {\it scheme} if $1_\Omega\in S$, and {\it commutative} if $c_{rs}^t=c_{sr}^t$ for all $r,s,t\in S$.  

\subsection{Isomorphisms and schurity} A  {\it combinatorial isomorphism} or, briefly, \emph{isomorphism} from $\cX$ to a coherent configuration  $\cX'=(\Omega', S')$ is defined to be a bijection $f: \Omega\rightarrow \Omega'$ such that  $s^f=\{(\alpha^f, \beta^f): (\alpha, \beta)\in s\}$ belongs to~$S'$ for all $s\in S$. In this case,  $\cX$ and $\cX'$ are said to be {\it isomorphic}; the set of all isomorphisms from $\cX$ to~$\cX'$ is denoted by $\iso(\cX,\cX')$. The group of all isomorphisms of $\cX$ to itself contains a normal subgroup
$$
\aut(\cX)=\{f\in\sym(\Omega):\ s^f=s \text{ for all } s\in S\}
$$
called the {\it automorphism group} of $\cX$.\medskip

Let $G\le\sym(\Omega)$. Denote by $S$ the set of all orbits $(\alpha,\beta)^G$ in the induced action of~$G$ on~$\Omega\times\Omega$, where $\alpha,\beta\in \Omega$. Then the pair $\inv(G)=(\Omega,S)$ is a coherent configuration. Any coherent configuration associated with a permutation group in this way is said to be {\it schurian}. Note that~$\cX$ is schurian if and only if $\cX=\inv(\aut(\cX)$.

\subsection{Extensions} There is a natural partial order\, $\le$\, on the set of all coherent configurations  on~$\Omega$. Namely, given two such coherent configurations~ $\cX$ and $\cY$, we set
$$
\cX\le\cY\ \Leftrightarrow\ S(\cX)^\cup\subseteq S(\cY)^\cup,
$$
and say that $\cY$ is the \emph{extension} of~$\cX$. The minimal and maximal elements with respect to this order are the trivial and {\it discrete} coherent configurations, respectively; in the last case, $S$ consists of singletons. Note that the functor $\cX\to\aut(\cX)$  reverse the inclusion, namely,
\qtnl{140322a}
\cX\le \cY \Rightarrow  \aut(\cX)\ge \aut(\cY).
\eqtn

\subsection{Algebraic isomorphisms and separability}\label{110622w} A bijection $\varphi:S\rightarrow S'$ is called an \emph{algebraic isomorphism} from $\cX$ to $\cX'$ if for all $r,s,t\in S$, we have
$$
c_{r^\varphi s^\varphi}^{t^\varphi}=c_{rs}^t.
$$
In this case, $|\Omega'|=|\Omega|$, $1^{}_{\Omega'}=\varphi(1_\Omega)$, and $\cX'$ is commutative if and only if so is $\cX$. Every $f\in\iso(\cX,\cX')$ induces the algebraic isomorphism $\varphi_f:\cX\to \cX',$ $s\mapsto s^f$; we put
$$
\iso(\cX,\cX',\varphi)=\{f\in\iso(\cX,\cX'):\ \varphi_f=\varphi\}.
$$
Note that $\aut(\cX)=\iso(\cX,\cX,\id)$, where $\id$ is the trivial (identical) algebraic automorphism of~$\cX$. 
Finally, for any $\alpha\in\Omega$ and $\alpha'\in\Omega'$, we put
$$
\iso_{\alpha,\alpha'}(\cX,\cX',\varphi)=\{f\in\iso(\cX_,\cX',\varphi):\ \alpha^f=\alpha'\}. 
$$

A coherent configuration $\cX$  is said to be \emph{separable} if every algebraic isomorphism from $\cX$  is induced by  isomorphism, equivalently, $\iso(\cX,\cX',\varphi)\ne\varnothing$ for all~$\cX'$ and~$\varphi$. Any trivial coherent configuration is separable.

\subsection{Relations}  The elements of $S$ and of $S^\cup$ are called {\it basis relations} and \emph{relations} of the coherent configuration~$\cX$. A unique basis relation containing the pair~$(\alpha,\beta)$ is denoted by~$r(\alpha,\beta)$.  The set of all relations is closed with respect to intersections and unions.

Any $\Delta\subseteq\Omega$ such that $1_\Delta\in S$ is called a {\it fiber} of $\cX$. In view of condition~(C1), the set~$F=F(\cX)$  of all fibers forms a partition of~$\Omega$. Every basis relation is contained in the Cartesian product of two uniquely determined fibers. When $\cX$ is a scheme, $F=\{\Omega\}$, whereas when $\cX$ is schurian, $F$ is just the set of orbits of the group~$\aut(\cX)$. 

Every algebraic isomorphism $\varphi:\cX\to\cX'$ is extended in a natural way to a poset isomorphism $S^\cup\to (S')^\cup$, denoted also by~$\varphi$. Then $\varphi$ induces 
a poset isomorphism $F^\cup\to (F')^\cup$, $\Delta\mapsto\Delta^\varphi$, where the set $\Delta^\varphi$ is defined by the equality $\varphi(1_\Delta)=1_{\Delta^\varphi}$.

Let $\cX\le \cY$ and $\cX'\le\cY'$. We say that the algebraic isomorphism $\varphi$ is extended to an algebraic isomorphism $\psi:\cY\to\cY'$  if $\psi(s)=\varphi(s)$ for all $s\in S$. Note that such an extension is uniquely defined if it exists. 

\subsection{Parabolics and quotients}\label{140322w}  A relation of $\cX$  that is an equivalence relation on $\Omega$ is called a \emph{parabolic} of~$\cX$; the set of all parabolics  is denoted by~$E=E(\cX)$. It contains the equivalence relations $\grp{s}$ and $\rad(s)$ for $s\in S^\cup$. Given a parabolic $e\in E$, we define $S_{\Omega/e}=\{s_{\Omega/e}:\ s\in S\}$, and $S_\Delta=\{s_\Delta:\ s\in S,\ s_\Delta\ne\varnothing\}$ for any $\Delta\in\Omega/e$. Then 
$$
\cX_{\Omega/e}=(\Omega/e,S_{\Omega/e})\qaq \cX_\Delta=(\Delta,S_\Delta)
$$  
are coherent configurations, called a \emph{quotient} of~$\cX$ modulo~$e$ and  a \emph{restriction} of~$\cX$ to~$\Delta$, respectively. Their fibers are obtained from fibers $\Gamma\in F$ as follows: $\Gamma_{\Omega/e}$ in the first case, and $\Gamma\cap\Delta$ in the second. In particular, $\cX_{\Omega/e}$ and $\cX_\Delta$ are schemes if so is~$\cX$.

Let $\varphi:\cX\to\cX'$ be an algebraic isomorphism and  $e\in E$. Then  $\varphi(e)$ is a parabolic of~$\cX'$ with the same number of classes. Moreover, 
\qtnl{030322a}
\varphi(\grp{s})=\grp{\varphi(s)}\qaq \varphi(\rad(s))=\rad(\varphi(s)).
\eqtn
The algebraic isomorphism $\varphi$ induces the algebraic isomorphism
$$
\varphi^{}_{\Omega^{}/e^{}}:\cX^{}_{\Omega^{}/e^{}}\to \cX'_{\Omega'/e'},\ s^{}_{\Omega^{}/e^{}}\mapsto s'_{\Omega'/e'},
$$
where $e'=\varphi(e)$ and $s'=\varphi(s)$.  Now let  $\Delta\in\Omega/e$ and $\Delta'\in\Omega/e'$. Assume that there is $\Gamma\in F$ such that $\Delta\cap\Gamma\ne\varnothing\ne \Delta'\cap\Gamma'$, where $\Gamma'=\Gamma^\varphi$. Then  $\varphi$ induces an algebraic isomorphism 
\qtnl{200322u}
\varphi^{}_{\Delta,\Delta'}:\cX^{}_{\Delta^{}}\to \cX'_{\Delta'},\ s^{}_{\Delta^{}}\mapsto s'_{\Delta'}.
\eqtn
This  algebraic isomorphisms always exists for all $\Delta$ and $\Delta'$ if $\cX$ (and hence $\cX'$) is a scheme.

\subsection{Sections} Let $\cX$ be a coherent configuration, $e\in E$, and $\Delta$ a class of a parabolic containing~$e$. The quotient set $\fS=\Delta/e$ is called a {\it section} of $\cX$. Any element of $\fS$ is of the form $\alpha_\fS=\alpha e$ for some $\alpha\in\Delta$; when $\Delta$ is implicit, we write $\alpha_\fS=\varnothing$  for all $\alpha\not\in \Delta$. For any $s\in S$, we define the relation $s_\fS$ on $\fS$ by formula~\eqref{130622a}.  The set of all sections of~$\cX$  is denoted by~$\fS(\cX)$. This set  is partially ordered: $\Delta/e\le \Delta'/e'$,  whenever $\Delta\subseteq\Delta'$ and $e'\subseteq e$.

\lmml{110322a}
Let $\cX$ be a commutative scheme, $\fS\in\fS(\cX)$, and $s\in S$. Then $\alpha_\fS s_\fS=(\alpha s)_{\fS}$ for all $\alpha\in\Omega$ such that $\alpha_\fS\ne\varnothing$.
\elmm
\prf Obviously,  $\alpha_\fS s_\fS\supseteq (\alpha s)_{\fS}$. Conversely, without loss of generality, we assume that $s_\fS\ne\varnothing$.  Let $\fS=\Delta/e$ and $\alpha\in\Delta$. Assume that $\beta_\fS\in \alpha_\fS s_\fS$ for some~$\beta\in\Delta$. Then there are $\alpha'\in\alpha e$ and $\beta'\in\beta e$ such that $(\alpha',\beta')\in s$. Then $c_{rs}^t\ne 0$ where $r=r(\alpha,\alpha')$ and $t=r(\alpha,\beta')$. By the commutativity, $c_{sr}^t\ne 0$. It follows that there exists $\beta''\in\beta e$ such that $(\alpha,\beta'')\in s$. Thus,  $\beta_\fS=\beta e=\beta''e\subseteq (\alpha s)_{\fS}$.\eprf

For a section $\fS=\Delta/e\in\fS(\cX)$, we put $\cX_\fS=(\cX_\Delta)_{\Delta/e}$. Then  $\cX_\fS$ is schurian if so is $\cX$.   Now let $\varphi:\cX\to\cX'$ be an algebraic isomorphism,  $e'=\varphi(e)$, and $\Delta'$  a class of the $\varphi$-image of the parabolic of~$\cX$, containing the class~$\Delta$.  The algebraic isomorphism~\eqref{200322u} (if it is defined)  induces the  algebraic isomorphism 
$$
\varphi_{\fS,\fS'}=(\varphi_{\Delta,\Delta'})_{\Delta/e,\Delta'/e'}
$$
between the coherent configurations  $\cX^{}_\fS$ and $\cX'_{\fS'}$.

\subsection{Point extensions}
The  {\it point extension} $\cX_{\alpha,\beta,\ldots}$ of  the coherent configuration~$\cX$ with respect to the points~$\alpha,\beta,\ldots\in\Omega$  is defined to be the smallest coherent configuration $\cY=(\Omega,T)$ such that $\cY\ge \cX$ and $1_\alpha,1_\beta,\ldots\in T$.  When the points are irrelevant, we use the term ``point extension''  and ``one-point extension'' if $\alpha=\beta=\ldots$.

Let $\varphi:\cX\to\cX'$ be an algebraic isomorphism, and let $\alpha\in\Omega^m$, $\alpha'\in{\Omega'}^m$. We say that algebraic isomorphism $\psi:\cX^{}_{\alpha^{}} \to\cX'_{\alpha'}$ extending $\varphi$ is  an {\it $(\alpha,\alpha')$-extension } of $\varphi$ if 
$$
\psi(1_{\alpha^{}_i})=1_{\alpha'_i},\quad i=1,\ldots,m.
$$
Note that the $(\alpha,\alpha')$-extension is unique if it exists.

When $m=1$,  the $(\alpha,\alpha')$-extension $\psi$ exists only if $\psi(r(\alpha,\alpha))=r(\alpha',\alpha')$, or, equivalently, $\Delta^\varphi=\Delta'$, where $\Delta$ (respectively,~$\Delta'$) is the fiber of $\cX$ (respectively,~$\cX'$), containing the point~$\alpha$ (respectively,~$\alpha'$).

 \subsection{Partly regular coherent configurations.}

A coherent configuration $\cX$ is said to be {\it partly regular} if there exists a point $\alpha \in \Omega$ such that $|\alpha s|\le 1$ for all $s \in S$; the point $\alpha$ is said to be~{\it regular}. In this case, the set $\alpha^{\aut(\cX)}$ is a faithful regular orbit of $\aut(\cX)$.  Note that every extension of partly regular  coherent configuration is partly regular.

\lmml{050222w1}{\rm\cite[Theorem 3.3.19]{CP2019}}
Every partly regular coherent configuration is schurian and separable.
\elmm

\crllrl{090322a}
Let $\cX$ be a partly regular coherent configuration and  $\cX'\ge \cX$. Then  $\cX'=\cX$ if and only if $F(\cX')=F(\cX)$.
\ecrllr
\prf
It suffices to verify the ``if'' part only. Assume that $F(\cX')=F(\cX)$. We take $\Delta\in F(\cX)$ that contains a regular point of~$\cX$. Since $\cX$ is schurian (Lemma~\ref{050222w1}), $\Delta$ is a faithful regular orbit of the group~$G=\aut(\cX)$ and hence $|G|=|\Delta|$. The same argument applied to partly regular coherent configuration $\cX'\ge\cX$ and $\Delta\in F(\cX')$, shows that $|G'|=|\Delta|$, where $G'=\aut(\cX')$.  Thus, 
$$
|G|=|\Delta|=|G'|.
$$
Since $G'\le G$ (see~\eqref{140322a}, this yields $G=G'$, and $\cX=\inv(G)=\inv(G')=\cX'$.
\eprf

\section{Multidimensional coherent configurations and  WL-dimension}\label{120622s}

The concept of the $\WL$-dimension of a graph was introduced in~\cite{Grohe2017} in terms of the multidimensional Weisfeiler-Leman algorithm. A goal of this section is to analyze this definition in terms of coherent configurations. Our approach is based on the multidimensional coherent configurations defined in~\cite{Babai2019}; all the necessary information and results about them are taken from~\cite{AndresHelfgott2017}.

Throughout this section, we fix an integer $m\ge 2$, and put $M=\{1,\ldots,m\}$. The monoid of all mappings $\sigma:M\to M$ is denoted by $\mon(M)$. Elements of the Cartesian $m$-power $\Omega^m$ are $m$-tuples $x=(x_1,\ldots,x_m)$ with $x_i\in\Omega$ for all $i\in M$. For $\alpha\in\Omega$, we put
\qtnl{030422q}
x_{i\leftarrow \alpha}=(x_1,\ldots,x_{i-1},\alpha,x_{i+1},\ldots,x_m).
\eqtn
For $X\subseteq\Omega^m$ and $\sigma\in M$,  we put $X^\sigma=\{x^\sigma:\ x\in X\}$, where $x^\sigma=(x_{1^\sigma},\ldots,x_{m^\sigma})$.

\subsection{The multidimensional Weisfeiler-Leman algorithm}
The $m$-dimensional Weisfeiler-Leman algorithm constructs for a given graph $\cG=(\Omega,D)$ a certain coloring $c(m,\cG)$ of the set $\Omega^m$; a \emph{coloring} is meant as a function from $\Omega^m$ to a linear ordered set the  elements of which are called \emph{colors}. At the first stage, an initial coloring $c_0=c_0(m,\cG)$ of $\Omega^m$ is defined by the following condition: given $x,y\in\Omega^m$, 
\qtnl{120622v}
c_0(x)=c_0(y)\quad\Leftrightarrow\quad
\{(i,j):\ (x_i,x_j)\in D\}=\{(i,j):\ (y_i,y_j)\in D\},
\eqtn
where the pair $(i,j)$ runs over the set $M\times M$. At the second stage, this coloring is refined step by step. Namely, if $c_i$ is the coloring constructed at the $i$th step ($i\ge 0$), then the color of an $m$-tuple $x$ in the coloring $c_{i+1}$ is defined to be 
$$
c_{i+1}(x)=(c_i(x),\mltst{(c_i(x_{1\leftarrow \alpha}),\ldots,c_i(x_{m\leftarrow \alpha})):\ \alpha\in\Omega}),
$$
where $\mltst{\cdot}$ denotes a multiset. The algorithm stops when $|\im(c_i)|=|\im(c_{i+1})|$ and the final coloring $c(m,\cG)$ is set to be~$c_i$.

The coloring $c_\cG=c(m,\cG)$ defines a partition $\WL_m(\cG)$ of $\Omega^m$ into the color classes $X=c_\cG^{-1}(i)$, where $i$ runs over the colors of~$c_\cG$. In particular,  the color $c_\cG(X)=c_\cG(x)$ does not depend on $x\in X$. Since the final coloring is a refinement of the initial one, we have the following statement.

\lmml{120622t}
Let $X\in\WL_m(\cG)$ and $i,j\in M$. Then $(x_i,x_j)\in D$  either for all or for no $x\in X$.
\elmm 

Another important property of the coloring $c_\cG$ is that it defines a monoid homomorphism $\tau:\mon(M)\to\mon(\im(c_\cG))$ such that for all $X\in\WL_m(\cG)$ and all $\sigma\in\mon(M)$, we have
\qtnl{120622t3}
c_\cG(X)^{\tau(\sigma)}=c_\cG(X^\sigma).
\eqtn
According with interpretation given in~\cite{AndresHelfgott2017}, this equality means that the color of~$X$ ``knows'' the color of $X^\sigma$.

Two graphs $\cG$ and $\cG'$ are said to be {\it $\WL_m$-equivalent}\footnote{In terms of \cite{Grohe2017}, this means that the $m$-dimensional Weisfeiler-Leman algorithm does not distinguish $\cG$ and~$\cG'$.} if $\im(c_{\cG^{}})=\im(c_{\cG'})$. The \emph{Weisfeiler-Leman dimension} $\dim_{\scriptscriptstyle\WL}(\cG)$ of a graph $\cG$ is defined to be the smallest natural $m$ such that every graph $\WL_m$-equivalent to~$\cG$ is isomorphic to~$\cG$. It should be noted that these definitions can be extended to the case $m=1$.

For $m=2$, the partition $\WL(\cG)=\WL_2(\cG)$ coincides with $S(\cX)$ for some coherent configuration on $\Omega$, called the coherent configuration of the graph~$\cG$; in fact, $\cX$ is the smallest coherent configuration for which $D\in S^\cup$. The graphs $\cG$ and $\cG'$ are $\WL_2$-equivalent (briefly, $\WL$-equivalent) if and only if the mapping $X\to c_{\cG'}^{-1}(c_{\cG^{}}(X))$ is an algebraic isomorphism from $\WL(\cG)$ to $\WL(\cG')$ (see~\cite{Fuhlbr2018a}), where $c_{\cG^{}}=c(2,\cG)$ and $c_{\cG'}=c(2,\cG')$. A goal of this section is to prove an analog of the  ``only if'' part of this statement for $m\ge 3$. 

\subsection{Multidimensional coherent configurations}\label{140422a}
For a tuple $x\in\Omega^m$, denote by $\rho(x)$ the equivalence relation on $M$ such that $(i,j)\in\rho(x)$ if and only if  $x_i=x_j$.   For  any $X_1,\ldots,X_m\subseteq \Omega^m$, denote by $n(x;X_1,\ldots,X_m)$ the number of all $\alpha\in\Omega$ such that $x_{i\leftarrow \alpha} \in X_i$ for all~$i\in M$.

\dfntn
A partition $\fX$ of  $\Omega^m$ is called an $m$-ary coherent configuration on~$\Omega$  if the following conditions are satisfied for all $X\in \fX$:
\nmrt
\tm{C1'} $\rho(X):=\rho(x)$ does not depend on $x\in X$,
\tm{C2'} $X^\sigma\in \fX$  for all $\sigma\in\mon(M)$,
\tm{C3'} for any $X_0,X_1,\ldots,X_m\in \fX$, the number $n_{X_1,\ldots,X_m}^{X_0}=n(x_0;X_1,\ldots,X_m)$ 
does not depend on $x_0\in X_0$.
\enmrt
\edfntn

For $m=2$, conditions (C1'),  (C2'), and (C3') imply conditions (C1), (C2), and (C3), respectively. In fact,  the coherent configurations are just the $2$-ary coherent configurations. An example of $m$-ary coherent configuration is a partition of $\Omega^m$ into the orbits of a permutation group on $\Omega$, acting on $\Omega^m$ coordinatewise.

The set of all $m$-ary configurations on $\Omega$ is partially ordered. Namely, $\fX\le \fY$ if every class of $\fX$ is a union of some classes of~$\fY$, or equivalently, if $\fX^\cup\subseteq \fY^\cup$, where $\fX^\cup$ (respectively, $\fY^\cup$) is the set of all unions of classes of~$\fX$ (respectively,~$\fY$). The largest $m$-ary coherent configuration is the discrete one in which every class is a singleton; the smallest $m$-ary coherent configuration consists of the orbits of the symmetric group $\sym(\Omega)$ in its componentwise action on~$\Omega^m$ (this easily follows from (C2) and the fact that the orbits are in one-to-one correspondence with the equivalence relations on~$M$).

An \emph{algebraic isomorphism} of $m$-ary coherent configurations $\fX$ and $\fX'$ is a bijection $\varphi:\fX\to\fX'$ such that for all $X,X_0,\ldots,X_m\in \fX$ and $\sigma\in\mon(M)$,
\qtnl{110622h}
\varphi(X^\sigma)=\varphi(X)^\sigma\qaq n_{X_1,\ldots,X_m}^{X_0}=n_{\varphi(X_1),\ldots,\varphi(X_m)}^{\varphi(X_0)}.
\eqtn
For $m=2$, one can easily verify that our definition is agreed with that in Subsection~\ref{110622w}.

Note that if   $X\in \fX$ and $i,j\in M$, then $(i,j)\in\rho(X)$ if and only if $X^\sigma=X$, where $\sigma\in\mon(M)$ is identical on $M\setminus\{i\}$ and takes $j$ to~$i$. It follows that
\qtnl{110622x}
\rho(\varphi(X))=\rho(X).
\eqtn
We extend $\varphi$ to a bijection $\fX^\cup\to (\fX')^\cup$ in a natural way; then $X\subseteq Y$ implies $\varphi(X)\subseteq \varphi(Y)$ for all $X,Y\in\fX^\cup$. 

An example of algebraic isomorphism of $m$-ary coherent configurations is obtained for  any two $\WL_m$-equivalent graphs~$\cG$ and $\cG'$. In this case, we have $\im(c_{\cG^{}})=\im(c_{\cG'})$ and  the mapping
\qtnl{110622l}
\varphi_{\cG,\cG'}:\WL_m(\cG)\to\WL_m(\cG'),\ X\mapsto c^{-1}_{\cG'}(c^{}_{\cG^{}}(X)),
\eqtn
is a bijection.  In view of equality~\eqref{120622t3}, it satisfies the first relation in~\eqref{110622h}, whereas the second relation  follows from the description of the $m$-dimensional Weisfeiler-Leman algorithm. Thus, $\varphi_{\cG,\cG'}$ is an algebraic isomorphism from
$\WL_m(\cG)$ to $\WL_m(\cG')$. 

\subsection{Projections}
Let $k\in M$. The \emph{$k$-projection} of $\Omega^m$ is defined to be the mapping 
$$
\pr_k:\Omega^m\to \Omega^k, \ (\alpha_1,\ldots,\alpha_m)\mapsto (\alpha_1,\ldots,\alpha_k).
$$
The statement below  follows from \cite[Exercises 2.7,2.11]{AndresHelfgott2017} (and their proofs). Below, for every $X\subseteq\Omega^k$, we put
\qtnl{110622w5}
\wh X=\{(x_1,\ldots,x_k,x_k,\ldots,x_k)\in \Omega^m:\ (x_1,\ldots,x_k)\in X\}.
\eqtn

\lmml{100622a}
Let $\fX$ be an $m$-ary coherent configuration. Then  
$$
\pr_k(\fX)=\{\pr_k(X):\ X\in\fX\}
$$ 
is a $k$-ary coherent configuration.\footnote{In~\cite{Babai2019,AndresHelfgott2017}, the partition $\pr_k(\fX)$ is called the $k$-skeleton of~$\fX$.}  Moreover,
\qtnl{110622j1}
n_{X_1,\ldots,X_m}^{X_0}=\sum_{Y_1,\ldots,Y_{m-k}\in\fX} c_{\wh X_1,\ldots,\wh X_k,Y_1,\ldots,Y_{m-k}}^{\wh X_0}
\eqtn
for all $X_0,X_1,\ldots,X_m\in\pr_k(\fX)$.
\elmm

Let $\cG=(\Omega,D)$ be a graph and $k=2$. From Lemma~\ref{120622t}, it follows that $\wh D\in \fX^\cup$. Therefore, $D=\pr_2(\wh D)$ is a relation of the coherent configuration $\pr_2(\fX)$. Thus,
\qtnl{110622d}
\pr_2(\WL_m(\cG))\ge \WL(\cG).
\eqtn

Let $\varphi:\fX\to\fX'$ be an algebraic isomorphism of $m$-ary coherent configurations. We define a mapping $\varphi_k:\pr_k(\fX)\to\pr_k(\fX')$ by setting
$$
\varphi_k(\pr_k(X))=\pr_k(\varphi(X))
$$ 
for all $X\in\fX$. Note that if $\pr_k(X)=\pr_k(Y)$ for some $Y\in\fX$, then $X^\sigma=Y^\sigma$, where $\sigma\in\mon(M)$ takes $i$ to $\min\{i,k\}$.  For this~$\sigma$, we have 
$$
\pr_k(\varphi(X))=\pr_k(\varphi(X)^\sigma)\qaq \pr_k(\varphi(Y))=\pr_k(\varphi(Y)^\sigma).
$$
Since also $\varphi(X)^\sigma=\varphi(X^\sigma)=\varphi(Y^\sigma)=\varphi(Y)^\sigma$, we conclude that $\pr_k(\varphi(X))=\pr_k(\varphi(Y))$. Thus the mapping $\varphi_k$ is well-defined. Reversing the above argument, one can see that it is injective. Finally, $\varphi_k$   is surjective, because it is a composition of the surjections~$\pr_k$ and~$\varphi$.  

\lmml{100622a1}
Let $\varphi:\fX\to\fX'$ be an algebraic isomorphism of $m$-ary coherent configurations. Then the bijection $\varphi_k:\pr_k(\fX)\to\pr_k(\fX')$ is an algebraic isomorphism.
\elmm
\prf
Let $X\in\pr_k(\fX)$ and $\sigma\in\mon(K)$ with $K=\{1,\ldots,k\}$. Then, obviously,
$$
X^\sigma=\pr_k(\wh X)^\sigma=\pr_k(\wh X^{\wh \sigma}),
$$ 
where  $\wh\sigma\in \mon(M)$ is the mapping identical on $M\setminus K$ and coinciding with $\sigma$ on~$K$. It follows that
$$
\varphi_k(X^\sigma)=\varphi_k(\pr_k(\wh X^{\wh \sigma}))=\pr_k(\varphi(\wh X^{\wh \sigma}))=
\pr_k(\varphi(\wh X)^{\wh \sigma})=(\pr_k(\varphi(\wh X))^\sigma=\varphi_k(X)^\sigma,
$$
which proves the first part of~\eqref{110622h}. To prove the second relation, let $X_i\in\pr_k(\fX)$, where $i=0$ or $i\in K$. In view of formula~\eqref{110622x}, we have $\rho(\varphi(\wh X_i))=\rho(\wh X_i)$. It follows that $\varphi(\wh X_i)=\wh{X'}$ for some $X'\in \pr_k(\fX')$. Moreover,
$$
X'=\pr_k(\wh X')=\pr_k(\varphi(\wh X_i))=\varphi_k(\pr_k(\wh X_i))=X'_i,
$$
where $X'_i=\varphi_k(X_i)$. Consequently, $\varphi(\wh X_i)=\wh{X'_i}$ for all~$i$. Now, the required statement follows from formula~\eqref{110622j1}.
\eprf

For arbitrary $K\subseteq M$ and $X\subseteq~\Omega^m$, one can define $\pr_K(X)=(\pr_k(X^\sigma))^{\sigma^{-1}}$, where $\sigma\in\sym(M)$ is such that $K^\sigma=\{1,\ldots,k\}$ with $k=|K|$. Using condition~(C2'), one can easily prove Lemmas~\ref{100622a} and~\ref{100622a1} with $\pr_k$ replaced by~$\pr_K$. 

\subsection{Residues}
Let again $k\in M$ and $y\in \Omega^{m-k}$. For arbitrary $X\subseteq \Omega^m$, we define   
the residue $X_y$ of $X$ with respect to~$y$ to be the set of all $k$-tuples $\bar x$ such that $(\bar x_1,\ldots,\bar x_k,y_1,\ldots,y_{m-k})\in X$. The statement below  follows from \cite[Exercise~2.13]{AndresHelfgott2017} (and its proof). 

\lmml{100622c}
Let $\fX$ be an $m$-ary coherent configuration. Then  for any $y\in \Omega^{m-k}$,
$$
\fX_y=\{X_y:\ X\in\fX,\ y\in\pr_k(X)\}
$$ 
is a $k$-ary coherent configuration. Moreover,
$\fX_y\ge \pr_k(\fX)$ and
\qtnl{120622a}
n_{X_1,\ldots,X_m}^{X_0}=\sum_{Y_1,\ldots,Y_{m-k}\in\fX} c_{X_1y,\ldots,X_ky,Y_1,\ldots,Y_{m-k}}^{X_0y}
\eqtn
for all $X_0,X_1,\ldots,X_m\in\fX_y$, where $X_iy\in\fX$  contains $(x_1,\ldots,x_k,y_1\ldots,y_{m-k})$ for some (and hence for all) $(x_1,\ldots,x_k)\in X_i$, $i=0,\ldots,k$.
\elmm

The $k$-ary coherent configuration $\fX_y$ contains the singleton $\{(y_i,\ldots,y_i)\}$, where  $i=1,\ldots,m-k$. For $k=2$, it is just $1_{y_i}$, and the coherent configuration $\WL_m(\cG)_y$ is larger than or equal to the smallest coherent configuration on $\Omega$, for which $D$ and  the $1_{y_i}$ are relations. The last coherent configuration is the point extension of $\WL(\cG)$ with respect to $y_1,\ldots,y_{m-2}$. Thus,
\qtnl{110622e}
(\WL_m(\cG))_y\ge \WL(\cG)_{y_1,\ldots,y_{m-2}}.
\eqtn

Let $\varphi:\fX\to\fX'$ be an algebraic isomorphism of $m$-ary coherent configurations. Assume that  $y'\in \varphi_{m-k}(Y)$, where $y\in Y\in\pr_{m-k}(\fX)$. Then by Lemma~\ref{100622a1}, we have for every $Z\in\fX$,
$$
\pr_{m-k}(Z)=Y\quad\Leftrightarrow\quad \pr_{m-k}(\varphi(Z))=Y'.
$$
In particular, the sets $\{Xy:\ X\in\fX_y\}$ and $\{X'y':\ X'\in\fX'_{y'}\}$ are in a one-to-one correspondence, where  $Xy$ and $X'y'$ are as in Lemma~\ref{100622c}.  We define a mapping 
$$
\varphi_{y,y'}:\fX^{}_{y{}}\to\fX'_{y'},\ X\mapsto \varphi(Xy)_{y'}.
$$
Note that for fixed $y$ and $y'$, the mappings $\fX_y\to\fX$, $X\mapsto Xy$, and $\fX'_{y'}\to\fX'$, $X'\mapsto X'y'$, are injective. By the choice of $y$ and $y'$, this shows that $\varphi_{y,y'}$ is a bijection.

\lmml{040622a}
Let $\varphi:\fX\to\fX'$ be an algebraic isomorphism of $m$-ary coherent configurations,   $Y\in\pr_{m-k}(\fX)$, and $Y'=\varphi_{m-k}(Y)$. Then for every $y\in Y$ and every $y'\in Y'$, the bijection $\varphi_{y,y'}:\fX^{}_{y^{}}\to\fX'_{y'}$ is an algebraic isomorphism. Moreover, it extends the algebraic isomorphism $\varphi_k$.
\elmm
\prf
Let $X\in\fX_y$ and $\sigma\in\mon(K)$ with $K=\{1,\ldots,k\}$. Then, obviously,
$X^\sigma y=(Xy)^{\wh \sigma}$, where  $\wh\sigma\in \mon(M)$ is the mapping identical on $M\setminus K$ and coinciding with $\sigma$ on~$K$. It follows that
$$
\varphi_{y,y'}(X^\sigma)=\varphi(X^\sigma y)_{y'}=\varphi((Xy)^{\wh \sigma})_{y'}=(\varphi(Xy)^{\wh\sigma})_{y'}=(\varphi(X_y)_{y'})^\sigma=\varphi_{y,y'}(X)^\sigma,
$$
which proves the first part of~\eqref{110622h}. To prove the second relation, it suffices to note that $\varphi(Xy)=\varphi_{y,y'}(X)_{y'}$ and make use of formula~\eqref{120622a}.
\eprf

\subsection{Reduction to coherent configurations}
The main result of this subsection (Theorem~\ref{080622a} below) establishes a necessary condition for two graphs to be $\WL_m$-equivalent in terms of their coherent configurations.

\thrml{080622a}
Let $\cG=(\Omega,D)$ and $\cG'=(\Omega',D')$  be $\WL_m$-equivalent graphs,  $m\ge 2$. Then there is an algebraic isomorphism $\varphi:\WL(\cG)\to\WL(\cG')$ such that $\varphi(D)=D'$.  Moreover, if  $x\in\Omega^m$ and $x'\in{\Omega'}^m$ are such that $c_{\cG^{}}(x)=c_{\cG'}(x')$, then $\varphi$ has the $(y,y')$-extension with $y=\pr_{m-2}(x)$ and $y'=\pr_{m-2}(x')$. 
\ethrm
\prf
Denote by $\psi$ the algebraic isomorphism defined by~\eqref{110622l}. By Lemma~\ref{100622a1} for $k=2$, there is an algebraic isomorphism $\psi_2:\pr_2(\WL_m(\cG))\to\pr_2(\WL_m(\cG'))$. Note that $D$ is a relation of the coherent configuration $\cX=\WL(\cG)$ (see formula~\eqref{110622d}), and also
$$
\psi_2(D)=\pr_2(\psi(\wh{D^{}})=\pr_2(\psi(\wh{D'})=D'.
$$
Since $\cX$ is the smallest coherent configuration containing $D$ as a relation, $\psi_2(\cX))$ is the smallest coherent configuration containing $D'$ as a relation, i.e.,  $\psi_2(\cX)$ is equal to $\cX'=\WL(\cG')$. Thus as the required algebraic isomorphism~$\varphi$ one can take the restriction of $\psi_2$ to $\cX$.

Now, let  $x\in\Omega^m$ and $x'\in{\Omega'}^m$ be such that $c_{\cG^{}}(x)=c_{\cG'}(x')$. Then the algebraic isomorphism~$\psi$ takes the class $X\ni x$ to the class $X'\ni x'$. It follows that 
$$
\psi_{m-2}(\pr_{m-2}(X))=\pr_{m-2}(X').
$$
Since $y\in \pr_{m-2}(X)$ and $y'\in \pr_{m-2}(X')$, we can apply Lemma~\ref{040622a} to find  an algebraic isomorphism 
$$
\psi_{y,y'}: \WL_m(\cG)_{y^{}}\to\WL_m(\cG')_{y'}
$$
which extends $\psi_2$ and hence extends $\varphi$. On the other hand, the inclusion~\eqref{110622e} yields 
$$
\WL_m(\cG)_{y^{}}\ge \cX_{y_1,\ldots,y_{m-2}}\qaq\WL_m(\cG')_{y'}\ge \cX'_{y'_1,\ldots,y'_{m-2}}.
$$
Furthermore, the pairs $(1,2),(1,i+2)$ belong to the equivalence relation $\rho(1_{y_i}y)$. By formula~\eqref{110622x}, this implies that these pairs also belong to $\rho(\psi(1_{y_i}y))$. It follows
that $\psi_{y,y'}(1_{y_i})=\psi(1_{y_i}y)_{y'}= 1_{y'_i}$ for $i=1,\ldots,m-2$. Thus, 
$$
\psi_{y,y'}(\cX^{}_{y^{}_1,\ldots,y^{}_{m-2}})=\cX'_{y'_1,\ldots,y'_{m-2}},
$$
and the restriction of $\psi_{y,y'}$ to $\cX_{y'_1,\ldots,y'_{m-2}}$ is the $(y,y')$-extension of~$\varphi$.
\eprf

\section{Sesquiclosed coherent configurations and algebraic isomorphisms}\label{160322f}

\subsection{Sesquiclosed algebraic isomorphisms} In this subsection,  we introduce a notion of a sesquiclosed algebraic isomorphism, which weakens the notion of an (algebraic) $2$-isomorphism studied in~\cite[Section~3.5.2]{CP2019}. 

\dfntnl{200322q}
An algebraic isomorphism $\varphi:\cX\to\cX'$ is said to be  sesquiclosed if $\varphi$ has  the $(\alpha,\alpha')$-extension for all  $\alpha\in\Delta\in F(\cX)$ and $\alpha'\in\Delta'\in F(\cX')$ such that $\Delta^\varphi=\Delta'$. 
\edfntn

Examples of sesquiclosed algebraic isomorphisms arise naturally in the context of $\WL_3$-equivalent graphs, namely the following statement is an almost immediate consequence of Theorem~\ref{080622a}.

\lmml{010622a} 
Let $\cG$ and $\cG'$ be $\WL_3$-equivalent  graphs. Assume that $\cG$ is vertex-transitive. Then there is a sesquiclosed   algebraic isomorphism $\varphi:\WL(\cG)\to\WL(\cG')$ such that $\varphi(D)=D'$, where $D$ and $D'$ are the arc sets of $\cG$ and $\cG'$, respectively.
\elmm
\prf By  the first part of Theorem~\ref{080622a} (for $m=3$), there is an  algebraic isomorphism $\varphi:\WL(\cG)\to\WL(\cG')$ such that $\varphi(D)=D'$.  It suffices to verify that $\varphi$ has $(\alpha,\alpha')$-extension for all $\alpha\in\Omega$ and $\alpha'\in\Omega'$. Because $\cG$ is vertex-transitive, $X=\diag(\Omega^3)$ is a class of~$\WL_3(\cG)$. Since $\cG$ and $\cG'$ be $\WL_3$-equivalent, $X'=\diag({\Omega'}^3)$ is a class of~$\WL_3(\cG')$, and $c_{\cG^{}}(X)=c_{\cG'}(X')$. It follows that $c_{\cG^{}}(x)=c_{\cG'}(x')$, where $x=(\alpha,\alpha,\alpha)$ and $x'=(\alpha',\alpha',\alpha')$. By the second part of Theorem~\ref{080622a},  this implies that $\varphi$ has $(\alpha,\alpha')$-extension.
\eprf

Not every algebraic isomorphism is sesquiclosed, e.g., a straightforward computation shows that a unique nontrivial algebraic automorphism of the antisymmetric scheme of rank~$3$ and degree~$15$ is not sesquiclosed. On the other hand, \cite[Lemma~3.5.25]{CP2019} shows that every algebraic $2$-isomorphism is sesquiclosed.

\lmml{200222a}
Let $\cX$ be a coherent configuration, $\varphi:\cX\to\cX'$ a sesquiclosed algebraic isomorphism, $\fS\in\fS(\cX)$, and $\fS'\in\fS(\cX')$. Assume that $\varphi$ induces an algebraic isomorphism $\psi:\cX^{}_{\fS^{}}\to \cX'_{\fS'}$. Then $\psi$  is sesquiclosed.
\elmm 
\prf 
Let $\bar\alpha\in\bar\Delta\in F(\cX^{}_{\fS^{}})$ and $\bar\alpha'\in\bar\Delta'\in F(\cX'_{\fS'})$ be such that $\bar\Delta^\psi=\bar\Delta'$. Then there exist $\Delta\in F(\cX)$ and $\Delta'\in F(\cX')$ for which 
$$
\Delta_\fS=\bar\Delta,\quad \Delta'_{\fS'}=\bar\Delta',\quad \Delta^\varphi=\Delta'.
$$
Let  $\alpha\in\Delta$ and $\alpha'\in\Delta'$ be such that $\alpha^{}_{\fS^{}}=\ov\alpha$ and $\alpha'_{\fS'}=\alpha'$. The algebraic isomorphism $\varphi$ being sesquiclosed has $(\alpha,\alpha')$-extension $\varphi_{\alpha,\alpha'}$. Since $\varphi$ induces $\psi$, one can define the algebraic isomorphism  
$$
\psi_{\ov\alpha,\ov\alpha'}=(\varphi_{\alpha,\alpha'})_{\fS,\fS'}
$$
from $(\cX^{}_{\alpha^{}})_{\fS^{}}$ to $(\cX'_{\alpha'})_{\fS'}$. For any $\ov s\in S(\cX_\fS)$, we have
$$
\psi_{\ov\alpha,\ov\alpha'}(\ov s)=
(\varphi_{\alpha,\alpha'})_{\fS,\fS'}(s_\fS)=
(\varphi_{\alpha,\alpha'}(s))_{\fS'}=
\varphi(s)_{\fS'}=\psi(s_\fS),
$$
where $s\in S$ is such that $s_\fS=\ov s$, and
$$
\psi_{\alpha,\alpha'}(1_{\bar\alpha})=(\varphi_{\alpha,\alpha'})_\fS((1_\alpha)_\fS)=(\varphi_{\alpha,\alpha'}(1_\alpha))_\fS=(1_{\alpha'})_\fS=1_{\bar\alpha'},
$$
Thus the algebraic isomorphism $\psi_{\bar\alpha,\bar\alpha'}$ is the $(\bar\alpha,\bar\alpha')$-extension  of the algebraic isomorphism~$\psi$.\eprf

A coherent configuration $\cX$ is said to be {\it sesquiseparable} if every sesquiclosed algebraic isomorphism is induced by isomorphism. Clearly, every separable coherent configuration is sesquiseparable; the converse is not true, see Subsection~\ref{280522c}. 

\lmml{280222k}
A coherent configuration $\cX$ is sesquiseparable if every one-point extension of  $\cX$ is separable.
\elmm
\prf 
Let $\varphi:\cX\to \cX'$ be a sesquiclosed algebraic isomorphism. Then it has a  $(\alpha,\alpha')$-extension $\varphi_{\alpha,\alpha'}$ for all appropriate $\alpha$ and $\alpha'$. Since the  coherent configuration~$\cX_\alpha$ is separable, $\varphi_{\alpha,\alpha'}$ is induced by a certain isomorphism~$f$. Then $s^f=\varphi_{\alpha,\alpha'}(s)=\varphi(s)$ for all $s\in S$. Hence, $f\in\iso(\cX,\cX',\varphi)$ and $\cX$ is sesquiseparable.\eprf

\subsection{Sesquiclosed coherent configurations} In parallel with sesquiclosed algebraic isomorphisms, we introduce a notion which weakens the notion of $2$-closed coherent configuration studied in~\cite[Section~3.5.3]{CP2019}. Namely, in Lemma~3.5.25 there, it was proved that a coherent configuration $\cX$ is $2$-closed only if for every $\alpha\in\Omega$ the following conditions are satisfied:
\nmrt
\tm{S1} $F(\cX_\alpha)=\{\alpha s:\ s\in S,\ \alpha s\ne\varnothing\}$,
\tm{S2} the trivial algebraic automorphism of~$\cX$ is sesquiclosed.
\enmrt

\dfntnl{210322w}
A coherent configuration $\cX$ satisfying {\rm (S1)} and {\rm (S2)} for all $\alpha$ is said to be  sesquiclosed.
\edfntn

Though conditions (S1) and (S2) are satisfied in the schurian case, not every sesquiclosed coherent configuration is schurian, for example, those are coherent configurations constructed in~\cite[Theorem~4.2.4]{CP2019}. The same example shows that the class of sesquiclosed coherent configurations is not invariant with respect to algebraic isomorphisms. 

\lmml{090322c}
Let $\cX$ be a sesquiclosed coherent configuration and $\fS\in\fS(\cX)$. Then the coherent configuration $\cX_\fS$ is sesquiclosed.
\elmm

\prf Condition (S2) is satisfied for $\cX=\cX_\fS$ by Lemma~\ref{200222a} with $\cX=\cX'$ and $\varphi=\id$. Assume that condition (S1) is not satisfied for $\cX_\fS$. Then there are $\alpha\in\Omega$ and $s\in S$ such that 
\qtnl{210322p}
\varnothing\ne\alpha_\fS s_\fS\not\in F((\cX_\fS)_{\alpha_\fS}).
\eqtn
Without loss of generality, we may assume that $\alpha s\ne\varnothing$. By condition~(S1), we have $\alpha s\in F(\cX_\alpha)$ and hence $(\alpha s)_\fS\in F((\cX_\alpha)_\fS)$. However, $(\cX_\alpha)_\fS\ge (\cX_\fS)_{\alpha_\fS}$. Therefore, $(\alpha s)_\fS$ is contained in some $\Delta\in F((\cX_\fS)_{\alpha_\fS})$. Since $(\alpha s)_\fS\subseteq \alpha_\fS s_\fS$, we have
\qtnl{210322n}
(\alpha s)_\fS\subseteq \Delta\subseteq  \alpha_\fS s_\fS.
\eqtn
In view of~\eqref{210322p}, there is $\Delta'\in  F((\cX_\fS)_{\alpha_\fS})$ such that $\Delta\ne\Delta'\subseteq \alpha_\fS s_\fS$. Take
any point $\alpha'$ for which $(\alpha' s)_\fS$ intersects $\Delta'$. Then as above
\qtnl{210322n1}
(\alpha' s)_\fS\subseteq \Delta'\subseteq \alpha_\fS s_\fS.
\eqtn
Note that $\alpha$ and $\alpha'$ lie in the same fiber of~$\cX$, because $\alpha s\ne\varnothing\ne\alpha's$ and $s\in S$. By condition (S2), the trivial algebraic automorphism of~$\cX$ has the $(\alpha,\alpha')$-extension; denote it by $\psi$. Then there is an  algebraic isomorphism $\psi_\fS:(\cX_{\alpha^{}})_\fS\to(\cX_{\alpha'})_\fS$ that induces the trivial algebraic automorphism of~$\cX_\fS$. Using formulas~\eqref{210322n} and~\eqref{210322n1}, we obtain
$$
\Delta^{\psi_\fS}\supseteq((\alpha s)_\fS)^{\psi_\fS}=((\alpha s)^\psi)_\fS=(\alpha' s)_\fS\subseteq\Delta',
$$
whence $\Delta=\Delta'$, a contradiction. 
\eprf

\lmml{010622g}
Let $\cX$ be a sesquiclosed coherent configuration and $\varphi:\cX\to\cX'$ a sesquiclosed algebraic isomorphism. Then the coherent configuration $\cX'$ is sesqui\-closed.
\elmm
\prf
Let $\alpha'\in\Omega'$ and  $s'\in S'$. The algebraic isomorphism $\varphi:\cX\to\cX'$ being sesquiclosed, has an $(\alpha,\alpha')$-extension $\psi$ for some $\alpha\in\Omega$ such that $\alpha s\ne\varnothing$, where $s=\varphi^{-1}(s')$.  Since the coherent configuration $\cX$ is sesquiclosed, $\alpha s\in F(\cX_\alpha)$. Consequently, the set $\alpha's'= \psi(\alpha s)$ belongs to $F(\cX'_{\alpha'})$. This shows that $\cX'$ satisfies condition~(S1).

Let $\alpha',\alpha''\in\Omega'$ lie in the same fiber of $\cX'$. The algebraic isomorphism~$\varphi$ being sesquiclosed, has an $(\alpha,\alpha')$-extension $\psi$ and $(\alpha,\alpha'')$-extension  $\psi'$ for some $\alpha\in\Omega$. It follows that $\psi^{-1}\psi'$ is the $(\alpha',\alpha'')$-extension of the trivial algebraic automorphism of~$\cX'$. This shows that $\cX'$ satisfies condition~(S2).
\eprf

When a coherent configuration $\cX$ is not sesquiclosed, one can construct a uniquely determined sesquiclosed extension $\cX'$ such that $\aut(\cX')=\aut(\cX)$.  However, the explicit definition of $\cX'$  is outside the scope of the present paper. 

\subsection{Partly regular sections} In this subsection, we prove two auxiliary lemmas for the sections of sesquiclosed schemes, the restriction to which is partly regular. They will be used in Subsection~\ref{010622w}.

\lmml{230222ay}
Let $\cX$ be a commutative sesquiclosed  scheme, $\fS_0,\fS\in\fS(\cX)$, and $\alpha$ a  point such that $\alpha_{\fS_0}\ne\varnothing$. Assume that $\fS_0\le\fS$ and $((\cX_\fS)_{\alpha_\fS})_{\fS_0}$ is partly regular. Then
\qtnl{300522a}
((\cX_\fS)_{\alpha_\fS})_{\fS_0}
=((\cX_\alpha)_\fS)_{\fS_0}.
\eqtn
\elmm

\prf The obvious inclusion $(\cX_\fS)_{\alpha_\fS}\le (\cX_\alpha)_\fS$ implies
$
((\cX_\fS)_{\alpha_\fS})_{\fS_0}\le((\cX_\alpha)_\fS)_{\fS_0}.
$ 
The coherent configuration on the left-hand side  is partly regular by the hypothesis. By Corollary~\ref{090322a}, it suffices to verify that every fiber $\Delta$ of $((\cX_\alpha)_\fS)_{\fS_0}$ is also a fiber of $((\cX_\fS)_{\alpha_\fS})_{\fS_0}$.

By condition (S1), one can find $s\in S$ such that  $\Delta=(\alpha s)_{\fS_0}$. Since $\cX$ is commutative,  Lemma~\ref{110322a} yields
$$
\Delta=(\alpha s)_{\fS_0}=((\alpha s)_{\fS})_{\fS_0} =(\alpha_\fS s_\fS)_{\fS_0}.
$$
On the other hand, the scheme $\cX_\fS$ is sesquiclosed by Lemma~\ref{090322c}. Therefore, the set $\alpha_\fS s_\fS$ is a fiber of  the coherent configuration $(\cX_\fS)_{\alpha_\fS}$. It follows that $(\alpha_\fS s_\fS)_{\fS_0}$  and hence $\Delta$ is a fiber of  the coherent configuration $((\cX_\fS)_{\alpha_\fS})_{\fS_0}$. \eprf

\lmml{050222j}
Let $\cX$ and $\cX'$ be schemes,   $\varphi:\cX\to\cX'$ an algebraic isomorphism having an $(\alpha,\alpha')$-extension  $\varphi_{\alpha,\alpha'}$ for some $\alpha$ and $\alpha'$, and $f\in \iso_{\alpha,\alpha'}(\cX,\cX',\varphi)$. Assume that 
\nmrt
\tm{i} $\cX$ is schurian,  
\tm{ii}$(\cX_\alpha)_\fS$  is partly regular for some $\fS\in\fS(\cX)$ such that $\alpha_\fS\ne\varnothing$. 
\enmrt
Then 
\qtnl{210322v}
\iso(\cX^{}_{\alpha^{}},\cX'_{\alpha'},\varphi^{}_{\alpha,\alpha'})^\fS
=
\iso((\cX^{}_{\alpha^{}})_{\fS^{}},(\cX'_{\alpha'})_{\fS'},(\varphi_{\alpha,\alpha'})_{\fS,\fS'}),
\eqtn
where $\fS'=\fS^f$.
\elmm
\prf From the hypothesis, it follows that $f^\fS$ belongs to  the sets on the left- and right hand sides of~\eqref{210322v}; in particular, these sets  are not empty. Consequently,
$$
|\iso(\cX^{}_{\alpha^{}},\cX'_{\alpha'},\varphi^{}_{\alpha,\alpha'})^\fS|
=|\aut(\cX^{}_{\alpha^{}})^\fS|,
$$
and
$$
|\iso((\cX^{}_{\alpha^{}})_{\fS^{}},(\cX'_{\alpha'})_{\fS'},(\varphi_{\alpha,\alpha'})_{\fS,\fS'})|=|\aut((\cX^{}_{\alpha^{}})_{\fS^{}})|
$$
and it suffices to verify that 
\qtnl{280222a}
\aut((\cX^{}_{\alpha^{}})_\fS)=\aut(\cX^{}_{\alpha^{}})^{\fS^{}}.
\eqtn
Obviously, $\aut(\cX^{}_{\alpha^{}})^{\fS^{}}\le \aut((\cX^{}_{\alpha^{}})_\fS)$. Conversely, by condition~(ii), the group $\aut((\cX^{}_{\alpha^{}})_\fS)$ has a faithful regular orbit~$\Delta$. Since the coherent configuration $(\cX_\alpha)_\fS$ is schurian (Lemma~\ref{050222w1}), $\Delta=\Lambda_\fS$ for some fiber~$\Lambda\in F(\cX_\alpha)$.  By the first part of condition~(i), $\Lambda$ is an orbit of the group $\aut(\cX_\alpha)$. It follows that $\Delta=\Lambda_\fS$ is an orbit of the group $\aut(\cX_\alpha)^\fS\le \aut((\cX^{}_{\alpha^{}})_\fS)$. This proves~\eqref{280222a}, because $\Delta$ is a faithful regular orbit of the last group.\eprf

\section{$e_1/e_0$-condition}\label{090222e}

\subsection{Definition} Let $\cX=(\Omega,S)$ be a scheme, and let $e_0$ and $e_1$ be parabolics of~$\cX$ such that $e_0\subseteq e_1$. Following \cite{EvdP2004a}, we say that  $\cX$ satisfies the \emph{$e_1/e_0$-condition}  if for all $s\in S$,
$$
s\cap e_1=\varnothing\quad\Rightarrow\quad e_0\subseteq\rad(s).
$$  
This is always true if $e_0=1_\Omega$ or $e_1=\bone_\Omega$; in the other cases, we say that the $e_1/e_0$-condition is satisfied {\it nontrivially}. When $e_0=e_1$, the condition exactly means that the scheme~$\cX$ is isomorphic to the wreath product $\cX_\Delta\wr\cX_{\Omega/e_0}$ for any $\Delta\in\Omega/e_0$. Finally, we note that if $e'_0\subseteq e_0$ and $e'_1\supseteq e_1$ are parabolics, then $\cX$ satisfies also the $e'_1/e'_0$-condition. 

\rmrkl{190322a}
The $e_1/e_0$-condition provides internal definition of the operation wedge product of association schemes, introduced and studied in~{\rm \cite{Muzychuk2009a}}. Namely,  $\cX$ satisfies the $e_1/e_0$-condition if and only if $\cX$ is isomorphic to the wedge product of the schemes $\cX_\Delta$, $\Delta\in \Omega/e_1$, and the scheme $\cX_{\Omega/e_0}$.
\ermrk

In the sequel, we simplify the notation as follows. First, we put $\Omega_0=\Omega/e_0$, $\cX_0=\cX_{\Omega/e_0}$, and $\Omega_1=\Omega/e_1$. Second, for  every $\Delta\in\Omega_1$, we put $\Delta_0=\Delta/e_0$. Finally, for an algebraic isomorphism $\varphi:\cX\to \cX'$, we put $\varphi_0=\varphi_{\Omega_0}$.

\subsection{Admissible pairs}  Let $\varphi:\cX\to \cX'$,  $s\mapsto s'$ be an algebraic isomorphism, where $\cX'=(\Omega',S')$. In view of formula~\eqref{030322a},  $\cX'$ satisfies  the $e'_1/e'_0$-condition, where $e'_0=\varphi(e_0)$ and $e'_1=\varphi(e_1)$. As above, we define the sets $\Omega'_0$, $\Omega'_1$, and $\Delta'_0$ for all $\Delta'\in\Omega'_1$. In what follows, we are interested in finding the set $\iso(\cX,\cX',\varphi)$; for $\varphi=\id$, this was done in~\cite{EvdP2004a}.

Suppose we are given a bijection $f_0:\Omega^{}_0\to\Omega'_0$  taking $(e^{}_1)_{\Omega^{}_0}$ to $(e'_1)_{\Omega'_0}$, and  bijections $f_\Delta:\Delta\to \Delta'$ for each $\Delta\in \Omega_1$, where $\Delta'$ is a uniquely determined class of $e'_1$, for which $(\Delta^{}_0)^{f_0}=\Delta'_0$.  The pair $P=(\{f_\Delta\}_{\Delta\in \Omega_1},f_0)$ is said to be \emph{$e_1/e_0$-admissible} if 
\qtnl{040322w}
(f_\Delta)^{\Delta_0}=(f_0)^{\Delta_0}
\eqtn
for all $\Delta\in \Omega_1$. In this case there exists a uniquely determined bijection $f:\Omega\to \Omega'$ for which  $f^{\Omega_0}=f_0$  and $f^\Delta=f_\Delta$ for all $\Delta\in \Omega_1$. We say that $f$  is induced by the pair~$P$.

\thrml{310122d}
Let  $\cX$ be a scheme satisfying the $e_1/e_0$-condition,   $\cX'$ a scheme on~$\Omega'$, and $\varphi:\cX\to \cX'$ an algebraic isomorphism.  Then $f\in\iso(\cX,\cX',\varphi)$ if and only if
$f$ is induced by an $e_1/e_0$-admissible pair $(\{f_\Delta\},f_0)$ such that
\qtnl{030322c}
f_0\in\iso(\cX^{}_0,\cX'_0,\varphi^{}_0)\qaq f_\Delta\in\iso(\cX^{}_{\Delta^{}},\cX'_{\Delta'},\varphi^{}_{\Delta,\Delta'})\ \,\text{for all}\ \Delta\in \Omega_1.
\eqtn
\ethrm
\prf The ``only if'' part  is clear, because every	$f\in\iso(\cX,\cX',\varphi)$  is obviously induced by the pair $(\{f^\Delta\},f^{\Omega_0})$ which satisfies the conditions~\eqref{030322c}. Conversely, suppose that $f$ is induced by an $e_1/e_0$-admissible pair $(\{f_\Delta\},f_0)$  satisfying conditions~\eqref{030322c}. We need to verify that $s^f=\varphi(s)$ for all $s\in S$. 

Let  $s\subseteq e_1$.  Then $s$ is equal to the union of $s_\Delta$, $\Delta\in\Omega_1$. By the second part of~\eqref{030322c}, this implies that
$$
s^f=
\bigcup_{\Delta\in \Omega_1} (s_\Delta)^{f_\Delta}=
\bigcup_{\Delta\in \Omega_1} \varphi_{\Delta,\Delta'}(s_\Delta)=
\bigcup_{\Delta'\in \Omega'_1} \varphi(s)_{\Delta'}=\varphi(s).
$$

Now, we may assume that $s\cap e_1=\varnothing$. Then  $e_0\subseteq\rad(s)$, because $\cX$ satisfies the $e_1/e_0$-condition. It follows that $s$ is equal to the union of all $\Gamma\times\Lambda$ for which $(\Gamma,\Lambda)$ lies in $s_0=s_{\Omega/e_0}$. Because $f^{\Omega/e_0}=f_0$, we obtain
$$
s^f 
= \bigl(\bigcup_{(\Gamma,\Lambda)\in s_0} \Gamma\times \Lambda\bigr)^f
= \bigcup_{(\Gamma^f,\Lambda^f)\in s_0^f} \Gamma^f\times \Lambda^f
=\bigcup_{(\Gamma',\Lambda')\in s'_0} \Gamma'\times \Lambda',
$$
where $s'_0$, $\Gamma'$, and $\Lambda'$ are the $f_0$-images of $s_0$, $\Gamma$, and $\Lambda$, respectively. As was noted above, $\cX'$ satisfies the $e'_1/e'_0$-condition with $e'_0=\varphi(e_0)$, and $e'_1=\varphi(e_1)$. Since also $f_0$ induces $\varphi_0$, the relation $\varphi(s)$ equals the union of all $\Gamma'\times\Lambda'$ with $(\Gamma',\Lambda')\in s'_0$. Thus, 
$$		
s^f=
\bigcup_{(\Gamma',\Lambda')\in s_0'} \Gamma'\times \Lambda'=
\varphi(s),
$$		
as required.\eprf

\subsection{General sufficient condition}  A quite natural sufficient condition for  the existence of admissible $e_1/e_0$-pair satisfying~\eqref{030322c} could be given by the equality 
$$
\iso(\cX^{}_0,\cX'_0,\varphi^{}_0)^{\Delta_0}=\iso(\cX^{}_{\Delta^{}},\cX'_{\Delta'},\varphi^{}_{\Delta,\Delta'})^{\Delta_0}
$$ 
for all $\Delta\in \Omega_1$ (cf., \eqref{040322w}). In this subsection, we weaken this equality by considering smaller sets of isomorphisms, while assuming that the algebraic isomorphism~$\varphi$ has a one-point extension.

\lmml{120222a}
Let  $\cX$ be a  scheme  on~$\Omega$, satisfying the $e_1/e_0$-condition,  $\cX'$ a scheme on $\Omega'$, and $\varphi:\cX\to \cX'$ an algebraic isomorphism.  Assume that for all $\alpha\in\Omega$, $\alpha'\in\Omega'$,
\qtnl{080222q}
\iso_{\alpha^{}_0,\alpha'_0}(\cX^{}_0,\cX'_0,\varphi^{}_0)^{\Delta_0}
=
\iso_{\alpha,\alpha'}(\cX^{}_1,\cX'_1,\varphi^{}_1)^{\Delta_0},
\eqtn
where  $\alpha_0=\alpha e_0$, $\alpha'_0=\alpha' e'_0$, and $\cX_1=\cX_\Delta$ with $\Delta\in\Omega_1$ containing~$\alpha$,  $\cX'_1=\cX'_{\Delta'}$ with $\Delta'\in\Omega'_1$ containing $\alpha'$,  and $\varphi^{}_1=\varphi_{\Delta,\Delta'}$. Then for all $\alpha\in\Omega$ and  $\alpha'\in\Omega'$,
\nmrt
\tm{i} $\iso_{\alpha,\alpha'}(\cX,\cX',\varphi)^{\Omega_0}= \iso_{\alpha^{}_0,\alpha'_0}(\cX^{}_0,\cX'_0,\varphi^{}_0)$, %
\tm{ii} $\iso_{\alpha,\alpha'}(\cX,\cX',\varphi)^{\Delta}= \iso_{\alpha,\alpha'}(\cX^{}_1,\cX'_1,\varphi^{}_1)$. 
\enmrt
\elmm
\prf In both cases, it suffices to verify only the inclusion $\supseteq$ under the assumption that the set on the right-hand side is nonempty. In the case~(i), we take arbitrary  $f_0\in \iso_{\alpha^{}_0,\alpha'_0}(\cX^{}_0,\cX'_0,\varphi^{}_0)$. For every  $\Delta\in\Omega_1$, we will define a certain bijection $f_\Delta:\Delta\to \Delta'$, where $\Delta'$ is a unique class of~$e'_1$, for which $(\Delta^{}_0)^{f_0}=\Delta'_0$. First, assume that $\alpha\in \Delta$. Then $\alpha'\in\Delta'$ and
$$
(f_0)^{\Delta_0}\in \iso_{\alpha_0,\alpha'_0}(\cX^{}_0,\cX'_0,\varphi^{}_0)^{\Delta_0}.
$$
In view of~\eqref{080222q}, one can find an isomorphism $f_\Delta\in\iso_{\alpha,\alpha'}(\cX^{}_1,\cX'_1,\varphi^{}_1)$ for which $(f_\Delta)^{\Delta_0}=(f_0)^{\Delta_0}$; in particular,  $\alpha^{f_\Delta}=\alpha'$.
Now if $\alpha\not\in\Delta$, then we define~$f_\Delta$ as above, but instead of $\alpha$ and $\alpha'$ we take arbitrary points $\delta\in\Delta$ and $\delta'\in(\delta_0)^{f_0}$, respectively.
 
The constructed family $\{f_\Delta:\ \Delta\in\Omega_1\}$ together  with
the isomorphism~$f_0$ forms an $e_1/e_0$-admissible pair. Let $f:\Omega\to\Omega'$ be the bijection induced by this pair.  By definition, $f$  satisfies conditions~\eqref{030322c} for all~$\Delta$. Hence, $f\in\iso(\cX,\cX',\varphi)$ by Theorem~\ref{310122d}. Since also $\alpha^f=\alpha'$ and $f^{\Omega_0}=f_0$,  we are done.

In the case~(ii), let $f_1\in\iso_{\alpha,\alpha}(\cX^{}_1,\cX'_1,\varphi^{}_1)$, where $\alpha\in\Delta$ and  $\alpha'\in \Delta'$.  In view of~\eqref{080222q}, one can find 
$$
f_0\in\iso_{\alpha^{}_0,\alpha'_0}(\cX^{}_0, \cX'_0,\varphi^{}_0)
$$ 
such that $(f_1)^{\Delta_0}=f_0$. Let $\Gamma\in\Omega_1$ be the class of~$e_1$, other that~$\Delta$. Denote by~$\Gamma'$ a unique class of~$e'_1$, for which $(\Gamma^{}_0)^{f_0}=\Gamma'_0$, and take arbitrary $\gamma\in\Gamma$. Then $(\gamma_0)^{f_0}=\gamma'_0$ for some $\gamma'_0\in \Gamma'$. Using ~\eqref{080222q} again, we find $f_\Gamma\in\iso_{\gamma,\gamma'}(\cX^{}_{\Gamma^{}}, \cX'_{\Gamma'},\varphi^{}_{\Gamma,\Gamma'})$ such that $(f_\Gamma)^{\Gamma_0}=(f_0)^{\Gamma_0}$.  

The constructed family $\{f_\Lambda:\ \Lambda\in\Omega_1\}$ together  with the isomorphism~$f_0$ forms an $e_1/e_0$-admissible pair. Let $f:\Omega\to\Omega'$ is the bijection induced by this pair.  By definition, $f$  satisfies conditions~\eqref{030322c}. Hence,  $f\in\iso(\cX,\cX',\varphi)$ by Theorem~\ref{310122d}. Since also $\alpha'=\alpha^{f_1}=\alpha^f$ and $f^\Delta=f_1$, we are done.\eprf

\crllrl{120222ae}
In the hypothesis and notation of Lemma~{\rm \ref{120222a}}, the following statements are equivalent for all $\alpha\in\Omega$ and $\alpha'\in\Omega'$:
\nmrt
\tm{1} $\iso_{\alpha,\alpha'}(\cX,\cX',\varphi)\ne\varnothing$,
\tm{2} $\iso_{\alpha^{}_0,\alpha'_0}(\cX^{}_0,\cX'_0,\varphi^{}_0)\ne\varnothing$,
\tm{3} $\iso_{\alpha,\alpha'}(\cX^{}_1, \cX'_1,\varphi^{}_1)\ne\varnothing$.
\enmrt
\ecrllr

\crllrl{a070222a}
In the hypothesis and notations of Lemma~{\rm \ref{120222a}}, assume that $\cX=\cX'$ and $\varphi=\id$. Then $\cX$ is schurian if and only if $\cX_0$ is schurian and $\cX_1$ is schurian for some~$\alpha$.
\ecrllr
\prf  It suffices to prove the ``if'' part only. Assume that $\cX_0$ is schurian and~$\cX_1$ is schurian for some $\alpha\in\Omega$. Since $\varphi=\id$, we have $\cX'_0=\cX^{}_0$ (but not necessarily $\cX'_1=\cX^{}_1$, because, in general, $\alpha\ne\alpha'$).  Let  $\alpha'\in\Omega$. By the schurity of $\cX_0$, there exists $f\in\aut(\cX_0)$ such that $(\alpha^{}_0)^f=\alpha'_0$, or, equivalently,
$$
\iso_{\alpha^{}_0,\alpha'_0}(\cX^{}_0,\cX^{}_0,\id)\ne\varnothing.
$$
By Lemma~\ref{120222a}(i), this yields $\iso_{\alpha^{},\alpha'}(\cX,\cX, \id)\ne\varnothing$. Since $\alpha'$ is arbitrary, this means that the group $\aut(\cX)=\iso(\cX,\cX,\id)$ is transitive. Moreover,  by statements~(i) and~(ii) of Lemma~\ref{120222a} for $\alpha=\alpha'$, we have, respectively,
$$
(\aut(\cX)_\alpha)^{\Omega_0}= \aut(\cX_0)_{\alpha_0}\qaq
(\aut(\cX)_\alpha)^{\Delta}= \aut(\cX_1)_\alpha.
$$
By the transitivity of $\aut(\cX)$ and the schurity of~$\cX_0$ and~$\cX_1$, this implies that $\cX$ is schurian.\eprf

\subsection{Concrete sufficient condition}\label{010622w} The key point in our general sufficiency condition established in Lemma~\ref{120222a} is equality~\eqref{080222q}. In general, checking this equality is not an easy task. The lemma below gives a more or less simple criterion to be used in Section~\ref{050322f}.

\lmml{240222a}
Let  $\cX$ be a  sesquiclosed commutative scheme  on~$\Omega$, satisfying the $e_1/e_0$-condition,  $\cX'$ a scheme on $\Omega'$, and $\varphi:\cX\to \cX'$ a sesquiclosed algebraic isomorphism.  In the notation of Lemma~{\rm\ref{120222a}}, assume that for $\alpha\in\Omega$ and $\alpha'\in\Omega'$,
\nmrt
\tm{i}  $\cX_0$ and $\cX_1$ are schurian, $\iso(\cX^{}_0,\cX'_0,\varphi_0)\ne\varnothing\ne\iso(\cX^{}_1,\cX'_1,\varphi_1)$,
\tm{ii} $((\cX_0)_{\alpha_0})_{\Delta_0}$ and $((\cX_1)_\alpha)_{\Delta_0}$ are partly regular.
\enmrt
Then  equality~\eqref{080222q} holds.
\elmm
\prf  
Let $f_0\in \iso(\cX^{}_0,\cX'_0,\varphi_0)$. Because $\cX_0$ is a schurian scheme, we may assume that $(\alpha_0)^{f_0}=\alpha'_0$. Then the algebraic isomorphism $\psi_0$ induced by~$f_0$ is 
the $(\alpha^{}_0,\alpha'_0)$-extension of the algebraic isomorphism~$\varphi_0$.  It easily  follows that 
\qtnl{230322a}
\iso_{\alpha_0,\alpha'_0}(\cX^{}_0,\cX'_0,\varphi^{}_0)^{\Delta_0}=
\iso((\cX^{}_0)_{\alpha^{}_0},(\cX'_0)_{\alpha'_0},\psi_0)^{\Delta_0}.
\eqtn
Next, by the first parts of conditions~(i) and~(ii), the scheme $\cX_0$ and algebraic isomorphism $\varphi_0$ satisfy the hypothesis of Lemma~\ref{050222j} for $\fS=\Delta^{}_0$, $\fS'=\Delta'_0$, and $\alpha=\alpha_0$. Therefore,
\qtnl{050322c}
\iso((\cX^{}_0)_{\alpha^{}_0},(\cX'_0)_{\alpha'_0},\psi_0)^{\Delta_0}=
\iso(((\cX^{}_0)_{\alpha^{}_0})_{\Delta_0},((\cX'_0)_{\alpha'_0})_{\Delta'_0},(\psi_0)_{\Delta^{}_0,\Delta'_0}).
\eqtn
Furthermore, the scheme $\cX$,  sections $\fS=\Omega_0$, $\fS_0=\Delta_0$, and point $\alpha$ satisfy the hypothesis of Lemma~\ref{230222ay}. Therefore,
\qtnl{050322a}
((\cX^{}_0)_{\alpha^{}_0})_{\Delta_0}=((\cX_\alpha)_{\Omega_0})_{\Delta_0} =(\cX_\alpha)_{\Delta_0}.
\eqtn 

By Lemma~\ref{010622g}, the scheme $\cX'$ is sesquiclosed. Furthermore, the isomorphism $f_0$  induces an isomorphism from $((\cX^{}_0)_{\alpha^{}_0})_{\Delta_0}$ to  the coherent configuration  $((\cX'_0)_{\alpha'_0})_{\Delta'_0}$, which is therefore partly regular. Thus,  the scheme $\cX'$,  sections $\fS=\Omega'_0$, $\fS_0=\Delta'_0$, and point $\alpha'$ satisfy the hypothesis of Lemma~\ref{230222ay}, whence
\qtnl{050322a7}
((\cX'_0)_{\alpha'_0})_{\Delta_0}=((\cX'_{\alpha'})_{\Omega'_0})_{\Delta'_0}. =(\cX'_{\alpha'})_{\Delta'_0}.
\eqtn 

The $(\alpha,\alpha')$-extension $\varphi_{\alpha,\alpha'}$ of the (sesquiclosed) algebraic isomorphism $\varphi$ takes $\alpha_0$ to $\alpha'_0$. Therefore, $\varphi_{\alpha,\alpha'}$ extends $\varphi_0$ and takes $\Delta^{}_0$ to $\Delta'_0$. It follows that the algebraic isomorphism  $\psi_{\Delta^{}_0,\Delta'_0}$ defined as the restriction of $\varphi_{\alpha,\alpha'}$ to $((\cX_0)_{\alpha_0})_{\Delta}$ coincides with $(\psi_0)_{\Delta^{}_0,\Delta'_0}$.  With taking 
equalities~\eqref{050322a} and~\eqref{050322a7} into account, we can continue formula~\eqref{050322c} as follows:
$$
\iso(((\cX^{}_0)_{\alpha^{}_0})_{\Delta_0},((\cX'_0)_{\alpha'_0})_{\Delta'_0},(\psi_0)_{\Delta^{}_0,\Delta'_0})
=
\iso((\cX^{}_{\alpha^{}})_{\Delta^{}_0},(\cX'_{\alpha'})_{\Delta'_0}, \psi_{\Delta^{}_0,\Delta'_0}).
$$
This together with~\eqref{230322a}, yields
$$
\iso_{\alpha_0,\alpha'_0}(\cX^{}_0,\cX'_0,\varphi^{}_0)^{\Delta_0}=
\iso((\cX^{}_{\alpha^{}})_{\Delta^{}_0},(\cX'_{\alpha'})_{\Delta'_0}, \psi_{\Delta^{}_0,\Delta'_0}).
$$
Arguing in a similar way with $\cX_0$ and $\varphi_0$ replaced by $\cX_1$ and $\varphi_1$, we obtain
$$
\iso_{\alpha,\alpha'}(\cX^{}_1,\cX'_1,\varphi^{}_1)^{\Delta_0}=
\iso((\cX^{}_{\alpha^{}})_{\Delta^{}_0},(\cX'_{\alpha'})_{\Delta'_0}, \psi_{\Delta^{}_0,\Delta'_0}).
$$
Comparing the expressions on the right-hand sides of the last two equalities, we get~\eqref{080222q}.\eprf

 \section{Circulant sesquiseparable  schemes}\label{050322f}
 
\subsection{Circulant schemes} In this subsection, we follow \cite[Sections~2.4 and~4.4]{CP2019}. Let $G$ be a cyclic group. A coherent configuration~$\cX$ on the elements of $G$ is said to be {\it circulant} if for every $g\in G$, the permutation 
$$
\rho_g:x\mapsto xg,\ x\in G,
$$
is an automorphism of~$\cX$; in particular, $\cX$ is a commutative scheme.  There is a one-to-one correspondence between circulant schemes and the Schur rings over cyclic groups (see, e.g.,~\cite{EvdP2003}). The reader should not be embarrassed that most of the results cited in the present paper were formulated and proved in the language of Schur rings. In the rest of this subsection, we will change our terminology slightly to match that used for circulant schemes (and Schur rings). In the sequel, $\cX$ is a circulant scheme on $G$, $S=S(\cX)$, $E=E(\cX)$, and $\alpha=1_G$.   

Every group $H\le G$ defines an equivalence relation $e=e(H)$ on $G$,  the classes of which are the cosets of $H$ in~$G$. In particular, $e(\alpha)=\diag(G\times G)$ and $e(G)=G\times G$. Note that $H=\alpha e$ and $H\le H'$ if and only if $e(H)\subseteq e(H')$. 

We say that  $H$ is an {\it $\cX$-group} if $e(H)\in E$. The mapping $H\mapsto e(H)$ defines  a one-to-one correspondence between the $\cX$-groups and parabolics of~$\cX$. The parabolics~$\grp{s}$ and $\rad(s)$, where $s\in S$, correspond to $e(H)$ with $H=\grp{\alpha s}$ and $H=\{h\in G:\ (\alpha s)h=\alpha s\}$, respectively. Note that the parabolic $\rad(s)=e(H)$, and hence the group $H$, does not depend on the relation $s=r(\alpha,g)$, where $g$ is a generator of~$G$.  This group is called the {\it radical} of $\cX$ and is denoted by~$\rad(\cX)$.

In the context of circulant schemes, we are interested in those  sections $\fS=\Delta/e$ for which $\alpha \in\Delta$. In this case, $\Delta=U$ and $e=e(L)$ for some $\cX$-groups $U$ and $L$ such that $L\le U$; in particular, $\fS=U/L$. The set of all such sections is denoted by~$\fS_0(\cX)$. It is important to note that every $\fS\in\fS_0(\cX)$ is treated as a cyclic group and the scheme $\cX_\fS$ is treated as a circulant scheme on~$\fS$. We say that a section $\fS$ is \emph{trivial} (respectively, \emph{primitive}) if the scheme $\cX_\fS$ is  primitive (respectively, contains no proper $\cX_\fS$-groups). It is known that every primitive section of composite order is trivial.

\subsection{Normal circulant schemes}\label{280522c}
The scheme $\cX$ is said to be {\it normal} if the group   $\grp{\rho_g:\ g\in G}$ is normal in $\aut(\cX)$, or equivalently, $\aut(\cX)_\alpha\le\aut(G)$. This condition is always satisfied when $G$ is of prime order and $\cX$ is not trivial. 

\lmml{010322d1}{\rm \cite[Theorem~6.1 and Corollary~6.2]{EvdP2003}}
Every normal circulant scheme~$\cX$ is schurian. Moreover, any one-point extension of~$\cX$ is partly regular.
\elmm

There are normal circulant schemes that are not separable, see, e.g.,  \cite[Section~4.5]{CP2019}. However, from Lemmas~\ref{010322d1} and~\ref{050222w1}, it follows that every normal circulant scheme is  schurian and any one-point extension of it is separable. By Lemma~\ref{280222k}, this proves the following statement.

\crllrl{280522a}
Every normal circulant scheme is sesquiseparable.
\ecrllr

The following auxiliary statement is crucial for the proof of Theorem~\ref{150222e} below, because it enables us to verify condition~(ii) in Lemma~\ref{240222a}.

\lmml{050222b}
Let $\cX$ be a  normal circulant scheme and $\alpha=1_G$. Then the coherent configuration $(\cX_\alpha)_\fS$  is partly regular for every $\fS\in\fS_0(\cX)$.
\elmm
\prf  Let $\fS=U/L$, where $U$ and $L$ are $\cX$-subgroups, and let $u$ be a generator of~$U$. It suffices to verify that $\bar u=uL$ is a regular point of the coherent configuration ~$(\cX_\alpha)_\fS$.  Assume on the contrary that  there exist  $x,x'\in U$ such that 
\qtnl{050222u}
\bar x\ne\bar x'\qaq \bar r(\bar x,\bar u)=\bar r(\bar x',\bar u),
\eqtn
where $\bar x=xL$,  $\bar x'=x' L$, and for all $\ov y,\ov z\in \fS$ we denote by $\bar r(\ov y,\ov z)$  the basis relation of the coherent configuration $(\cX_\alpha)_\fS$, containing the pair $(\ov y,\ov z)$.  By the second equality in~\eqref{050222u}, there is $u'\in \bar u$ such that
$$
r(x,u)=r(x',u').
$$ 

The normality of $\cX$ implies that every element of $\aut(\cX_\alpha)=\aut(\cX)_\alpha$ is induced by raising to a power coprime to the order $n$ of the underlying cyclic group~$G$. Furthermore, the coherent configuration $\cX_\alpha$ is partly regular by Lemma~\ref{010322d1} and hence is schurian by Lemma~\ref{050222w1}. Thus there is an integer~$m$ coprime to~$n$ and such that 
$$
x^m=x'\qaq u^m=u'.
$$
In particular, $\bar x^m=\bar x'$ and $\bar u^m=\bar u'=\bar u$. On the other hand, since $u$ is a generator of~$U\ge L$,  there is an integer~$a$ for which $x=u^a$. Furthermore, $(uL)^m=\bar u^m=\bar u=uL$ and hence $u^m=yu$ for some $y\in L$. Thus for every $\ell\in L$, we have
$$
(x\ell)^m=(u^a \ell)^m=(u^a)^m  \ell^m=  (u^m)^a \ell^m=(yu)^a \ell^m=u^a(y^a \ell^m)=x\ell',
$$
where the element $\ell'=y^a \ell^m$ belongs to $L$. Consequently, $\bar x^m=\bar x$, implying  $\bar x'=\bar x^m=\bar x$ in contrast to the first inequality in~\eqref{050222u}.\eprf

Let $\cX$ be an arbitrary circulant scheme. A section $\fS\in\fS_0(\cX)$ is said to be \emph{normal} if the (circulant) scheme $\cX_\fS$ is normal, and \emph{subnormal} if there exists a normal section $\fS'\in\fS_0(\cX)$ such that $\fS\le\fS'$. Note that every normal section is subnormal and any primitive section of composite order greater than~$4$ is not subnormal.

\subsection{The $\fS$-wreath product of circulant schemes} Let $\cX$ be a circulant scheme and $\fS\in\fS_0(\cX)$. Assume that $\cX$ satisfies the $e(U)/e(L)$-condition, where~$U$ and~$L$ are $\cX$-groups such that $\fS=U/L$. In this case, we say that $\cX$ is the $\fS$-wreath product of schemes $\cX_1=\cX_U$ and $\cX_0=\cX_{G/L}$. The $\fS$-wreath product is nontrivial if $1<L\le U<G$, i.e., if  $\cX$ satisfies the $e(U)/e(L)$-condition nontrivially.  When $U=L$, we say that $\cX$ is a wreath product.

The main result of this subsection is to establish a simple criterion for the $\fS$-wreath product  to be  schurian and sesquiseparable. 

\thrml{150222e}
Let $\cX$ be a   circulant sesquiclosed scheme on $G$.  Assume that $\cX$ is the $\fS$-wreath product of (circulant) schemes $\cX_0$ and $\cX_1$ for some $\fS\in\fS_0$, and also 
\nmrt
\tm{1} the section  $\fS$ is subnormal in both  $\cX_0$ and $\cX_1$,
\tm{2} the schemes $\cX_0$ and $\cX_1$ are schurian  and sesquiseparable.
\enmrt
Then $\cX$ is schurian and sesquiseparable. 
\ethrm
\prf 
By condition~(1), the section $\fS$ is contained in certain  sections $\fS_0\in \fS(\cX_0)$ and $\fS_1\in \fS(\cX_1)$ such that the schemes
$$
\cX_{\fS_0}=(\cX_0)_{\fS_0}\qaq\cX_{\fS_1}=(\cX_1)_{\fS_1}
$$ 
are normal circulant.  By Lemma~\ref{010322d1},  the one-point extensions  $(\cX_{\fS_0})_\alpha$ and $(\cX_{\fS_1})_\beta$ are partly regular, where $\alpha=1_{\fS_0}$ and  $\beta=1_{\fS_1}$ treated as points of the schemes $\cX_{\fS_0}$ and $\cX_{\fS_1}$, respectively. 
By Lemma~\ref{050222b}, the coherent configurations 
$$
 ((\cX_{\fS_0})_\alpha)_\fS\le ((\cX_0)_\alpha)_\fS\qaq
((\cX_{\fS_1})_\beta)_\fS\le ((\cX_1)_\beta)_\fS
$$
are partly regular. It follows that so are the coherent configurations  $((\cX_0)_\alpha)_\fS$ and $((\cX_{\fS_1})_\beta)_\fS$. Consequently, condition~(ii) of Lem\-ma~\ref{240222a} is satisfied for $e_0=e(L)$, $e_1=e(U)$, and $\alpha=1_G$. Since $\aut(\cX)$ is transitive, this condition is true for all points of $\cX$. 

The schemes $\cX_0$ and $\cX_1$ are schurian by the first part of condition~(2). Thus condition~(i) of Lemma~\ref{240222a} is satisfied for $\varphi=\id$ (the condition $\iso(\cX^{}_1,\cX'_1,\varphi_1)\ne\varnothing$ is true, because the group $\aut(\cX)$ is transitive). By that lemma, equality~\eqref{080222q} holds for all $\alpha$ and $\alpha'$, and hence  the scheme $\cX$ satisfies the hypothesis of Corollary~\ref{a070222a}. Thus, $\cX$ is schurian.

Let $\varphi:\cX\to\cX'$ be a sesquiclosed algebraic isomorphism. Then the algebraic isomorphisms $\varphi_0$ and $\varphi_1$ are also sesquiclosed by Lemma~\ref{200222a}. By the second part of condition~(2), this implies that $\iso(\cX^{}_0,\cX'_0,\varphi_0)\ne\varnothing$, and $\iso(\cX^{}_1,\cX'_1,\varphi_1)\ne\varnothing$ . Thus condition~(i) of Lemma~\ref{240222a} is satisfied. By that lemma, equality~\eqref{080222q} holds for all $\alpha\in\Omega$ and $\alpha'\in\Omega$. Thus, $\iso(\cX,\cX',\varphi)\ne\varnothing$ by Corollary~\ref{120222ae}, and hence $\cX$ is sesquiseparable.\eprf

The conditions of Theorem~\ref{150222e} are always true when the circulant schemes $\cX_0$ and $\cX_1$ are normal. Thus this theorem implies the following statement.

\crllrl{190322c}
Let $\cX$ be a   circulant sesquiclosed scheme on $G$.  Assume that $\cX$ is the $\fS$-wreath product for some $\fS\in\fS_0(\cX)$  such that the both operands are normal. Then $\cX$ is schurian and sesquiseparable. 
\ecrllr

\subsection{Circulant schemes of prime power degree}  In this subsection, $\cX$ is a circulant scheme on a cyclic $p$-group~$G$. The structure of $\cX$  has been completely described in~\cite{Poschel1974} and~\cite{GolKN1981}  for $p\ge 3$ and $p=2$, respectively. In particular,  $\cX$ is schurian. 
Another important consequence of the description is the radical monotonicity for $p\ge 3$, namely,
\qtnl{250522a1}
\grp{s}\le\grp{t}\quad \Rightarrow\quad \rad(s)\le\rad(t).
\eqtn
for all $s,t\in S$. Using this fact, we show in the lemma below that the wreath product can be lifted from any section of $\cX$. 

\lmml{250522a2} 
Let $p\ge 3$ and $\fS\in\fS_0(\cX)$. If $\cX_\fS$  is  a nontrivial wreath product, then  so is $\cX$.
\elmm
\prf
Let $\fS=U/L$ and $\cX_\fS$  a nontrivial wreath product. Then there is an $\cX$-subgroup $H$ such that $L<H<U$ and for every $s\in S$,
$$
s\subseteq e(U)\setminus e(H)\quad\Rightarrow\quad e(H/L)\subseteq \rad(s_\fS),
$$
Without loss of generality, we may assume that $H$ is the only $\cX$-subgroup strictly between $L$ and $U$ (here we use the fact that $G$ is a cyclic $p$-group and hence the $\cX$-groups are linearly ordered with respect to inclusion).  It suffices to verify that
\qtnl{260522w}
s\subseteq e(U)\setminus e(H)\quad\Rightarrow\quad \rad(s)\ge e(H).
\eqtn  
Indeed, then by formula~\eqref{250522a1}, we have $\rad(s)\ge e(H)$ for all $s\in S$ not  contained in $e_H$. But this means that~$\cX$  satisfies the $e(H)/e(H)$-condition, i.e., $\cX$ is a nontrivial wreath product.

Suppose on the contrary  that \eqref{260522w} is not true for some  $s$. Then the $\cX$-subgroup~$L'$,  for which $\rad(s)=e(L')$, is strictly contained in $H$. By the assumption on~$H$, this implies that $L'<L$, and hence $\fS$  is a subsection of $\fS'=U/L'$. Note that 
$$
\rad(\cX_\fS)\ne 1\qaq \rad(\cX_{\fS'})=1.
$$
However, this contradicts the fact (see \cite[Corollary~7.4]{EvdokimovKP2013}) that if a circulant scheme (in our case $\cX_{\fS'}$) of odd prime power degree is of trivial radical, then every its section (in our case $\fS$) is also of trivial radical.
\eprf

In the following lemma, we collect some known properties of circulant schemes of prime power degree, that will be used in the proof of Theorem~\ref{010322a}. Below, we say that a circulant scheme $\cX$ is \emph{dense} if $\fS_0(\cX)$  contains no trivial section of composite order.

\lmml{250522a} 
 Let   $\cX$ be a circulant scheme of prime power degree. Then 
\nmrt
\tm{1}  if $\rad(\cX)=1$, then $\cX$ is trivial or normal,
\tm{2}  if  $\rad(\cX)\ne 1$ and $\cX$ is dense, then $\cX$ is a nontrivial $\fS$-wreath product,  and also $\rad(\cX_\fS)=1$ or $|\fS|=4$,
\tm{3} if $\fS_0(\cX)$ contains a proper primitive non-subnormal section, then $\cX$ is  a nontrivial wreath product.
\enmrt
\elmm   
\prf 
Statements~(1) and (3) are consequences of  \cite[Corollary~6.4]{EvdP2003},  and the statement cited  in~\cite[p.29]{Evdokimov2015} (and proved in \cite[Theorem 4.6]{EvdP2012a}), specified for schemes of prime power degree. Statement~(2) is exactly \cite[Lemma~6.1]{Evdokimov2015}.
\eprf
 
Now we are ready to prove the main result of this subsection.
 
\thrml{010322a}
Every circulant scheme of prime power degree is sesquiseparable.
\ethrm
\prf Let $\cX$ be a circulant scheme of degree $p^n$, where $p$ is a prime and $n\ge 1$. Assume first that  $\cX$ is trivial or normal. In the former case, $\cX$ is separable and hence sesquiseparable. In the latter case,  $\cX$ is sesquiseparable by Corollary~\ref{280522a}.\medskip
 
{\bf Claim.} {\it Let $\cX$ be neither trivial nor normal.  Then $\cX$ is a nontrivial $\fS$-wreath product with $\fS$ subnormal in both operands.}

\prf
By Lemma~\ref{250522a}(3), we may assume that  $\fS_0(\cX)$ contains no proper primitive non-subnormal section (otherwise we have a section $\fS$ for which $|\fS|=1$).  Then
$\fS_0(\cX)$ contains no trivial section of composite order, i.e., $\cX$ is dense. By statements~(1) and~(2) of Lemma~\ref{250522a}, the scheme $\cX$ is a nontrivial $\fS$-wreath product, and also $\rad(\cX_\fS)=1$ or $|\fS|=4$. In the last case, $\fS$ is normal (because all circulant schemes of degree at most~$4$ are normal) and hence is subnormal in both operands.

Let $\rad(\cX_\fS)=1$. Then by Lemma~\ref{250522a}(1), the section $\fS$ is either normal or trivial. In the first case, we are done as above. In the second case, $|\fS|=p\ge 5$ (otherwise $\fS$ is normal). It remains to verify that $\fS$ is subnormal in both operands. But if this is not true, then one of the operands is a nontrivial wreath product by Lemma~\ref{250522a}(3). Then so is $\cX$ by Lemma~\ref{250522a2}, and we are done as above. The claim is proved.
\eprf

To complete the proof of the theorem, we may assume by claim that $\cX$ is a nontrivial $\fS$-wreath product of circulant schemes~$\cX_0$ and $\cX_1$ such that the section~$\fS$ is subnormal both in $\cX_0$ and $\cX_1$. Note that $\cX_0$ and $\cX_1$ are circulant schemes of prime power degrees, which are less than $p^n$. By induction, this implies that they are sesquiseparable. Because they are also schurian, the scheme~$\cX$ is  sesquiseparable by Theorem~\ref{150222e}.\eprf

\subsection{Proof of Theorem~\ref{280222j}}
Let $\cG$ be a circulant graph of prime power order; in particular, $\cG$ is vertex-transitive and $\WL(\cG)$ is a schurian circulant scheme. Let $\cG'$  be a graph $\WL_3$-equivalent to~$\cG$. By Lemma~\ref{010622a}, there is a sesquiclosed   algebraic isomorphism $\varphi:\WL(\cG)\to\WL(\cG')$  such that $\varphi(D)=D'$, where $D$ and $D'$ are the arc sets of $\cG$ and $\cG'$, respectively. On the other hand, the scheme $\WL(\cG)$ is  sesquiseparable by Theorem~\ref{010322a}. Thus the algebraic isomorphism $\varphi$ is induced by an isomorphism~$f$. It follows that  $D^f=\varphi(D)=D'$. Hence,  $f\in\iso(\cG,\cG')$, i.e., the graphs $\cG$ and $\cG'$ are isomorphic. Consequently,   $\dim_{\scriptscriptstyle\WL}(\cG)\le 3$.

\providecommand{\bysame}{\leavevmode\hbox to3em{\hrulefill}\thinspace}
\providecommand{\MR}{\relax\ifhmode\unskip\space\fi MR }
\providecommand{\MRhref}[2]{%
	\href{http://www.ams.org/mathscinet-getitem?mr=#1}{#2}
}
\providecommand{\href}[2]{#2}

\end{document}